\documentstyle{amlts}
\begin{document}
\annalsline{157}{2003}
\received{May 18, 2001}
\startingpage{575}
\def\bye{\end{document}}
 \font\tenrm=cmr10
\def\joinrel{\mathrel{\mkern-4mu}}
\def\relbar{\mathrel{\smash-}}
\def\lrar{\relbar\joinrel\relbar\joinrel\relbar\joinrel\relbar\joinrel\relbar\joinrel\rightarrow}
\input boxedeps.tex 
\SetepsfEPSFSpecial 
\HideDisplacementBoxes
\def\figin#1#2{
$$
 {\BoxedEPSF{#1.eps scaled
#2}%
}%
$$
\noindent}

\catcode`\@=11
\font\twelvemsb=msbm10 scaled 1100
\font\tenmsb=msbm10
\font\ninemsb=msbm10 scaled 800
\newfam\msbfam
\textfont\msbfam=\twelvemsb  \scriptfont\msbfam=\ninemsb
  \scriptscriptfont\msbfam=\ninemsb
\def\msb@{\hexnumber@\msbfam}
\def\Bbb{\relax\ifmmode\let\next\Bbb@\else
 \def\next{\errmessage{Use \string\Bbb\space only in math
mode}}\fi\next}
\def\Bbb@#1{{\Bbb@@{#1}}}
\def\Bbb@@#1{\fam\msbfam#1}
\catcode`\@=12

 \catcode`\@=11
\font\twelveeuf=eufm10 scaled 1100
\font\teneuf=eufm10
\font\nineeuf=eufm7 scaled 1100
\newfam\euffam
\textfont\euffam=\twelveeuf  \scriptfont\euffam=\teneuf
  \scriptscriptfont\euffam=\nineeuf
\def\euf@{\hexnumber@\euffam}
\def\frak{\relax\ifmmode\let\next\frak@\else
 \def\next{\errmessage{Use \string\frak\space only in math
mode}}\fi\next}
\def\frak@#1{{\frak@@{#1}}}
\def\frak@@#1{\fam\euffam#1}
\catcode`\@=12


\newcommand{\Rr}{\Bbb R}
\newcommand{\Zz}{\Bbb Z}
\newcommand{\norm}[1]{\left\Vert#1\right\Vert}
\newcommand{\abs}[1]{\left\vert#1\right\vert}
\newcommand{\set}[1]{\left\{#1\right\}}
\newcommand{\eval}[1]{\left\langle#1\right\rangle}
\newcommand{\eps}{\varepsilon}
\newcommand{\To}{\longrightarrow}
\newcommand{\rmap}{\longrightarrow}
\newcommand{\Boxe}{\raisebox{.8ex}{\framebox}}
\newcommand{\X}{{{\cal X}}}
\newcommand{\F}{{{\cal F}}}
\newcommand{\U}{{{\cal U}}}
\newcommand{\Pa}{{{\cal P}}}
\newcommand{\NN}{{{\cal N}}}

\newcommand{\G}{{\cal G}}            
\newcommand{\s}{{\bf s}}             
\renewcommand{\t}{{\bf t}}           
\renewcommand{\H}{{\cal H}}          
\newcommand{\Mon}{{\rm Mon}\,}     
\newcommand{\Rep}{{\rm Rep}\,}     
\newcommand{\Hol}{{\rm Hol}\,}     

\newcommand{\A}{A}                      
\newcommand{\al}{\alpha}                
\newcommand{\be}{\beta}                 
\newcommand{\ga}{\gamma}                
\newcommand{\Lie}{{\cal L}}          
\renewcommand{\gg}{{\frak g}}        
\newcommand{\Exp}{{\rm Exp}\,}     
\newcommand{\Ker}{{\rm Ker}\,}     
\newcommand{\Ad}{{\rm Ad}\,}       
\newcommand{\ad}{{\rm ad}\,}       
\newcommand{\tr}{{\rm tr}\,}       

\renewcommand{\L}{\Bbb L}              

\newcommand{\flowA}[1]{\phi_{#1}^{s,t}}

\newcommand{\flowM}[1]{\Phi_{#1}^{s,t}} 
\newcommand{\flowG}[1]{\varphi_{#1}^{s,t}}

\newcommand{\rank}{{\rm rank}\,}   
\newcommand{\codim}{{\rm codim}\,}  


\title{Integrability of Lie brackets}
\shorttitle{Integrability of Lie brackets}
 \acknowledgements{The first author was supported in part by NWO and a Miller
Research Fellowship. The second author was supported in part by FCT
through program POCTI and grant POCTI/1999/MAT/33081.\hfill\break
\hglue22pt {\it Key words and phrases}.
 Lie algebroid, Lie groupoid.}
 \twoauthors{Marius Crainic}{Rui Loja Fernandes}
  
\institutions{University of California,
Berkeley, CA and Utrecht University,  Utrecht,\hfill\break The
Netherlands\\
{\eightpoint {\it E-mail address\/}: crainic@math.uu.nl}\\ \vskip6pt
Instituto Superior T\'{e}cnico,  Lisboa, Portugal\\
{\eightpoint {\it E-mail address\/}:  rfern@math.ist.utl.pt}}

\centerline {\bf Abstract}
\vglue12pt  In this paper we present the solution to a longstanding
problem of differential geometry: Lie's third theorem for Lie algebroids.
We show that the integrability problem is controlled by two
computable obstructions. As
applications we derive, explain and improve the known integrability
results, we establish integrability by local Lie groupoids, we clarify
the smoothness of the Poisson sigma-model for Poisson manifolds, and
we describe other geometrical applications. 
 
\def\sni#1{\smallbreak\noindent{#1}. }
\def\ssni#1{\vglue-1pt\noindent\hskip18pt {#1}.}

\vglue12pt \centerline{\bf Contents}

\sni{0} Introduction
\sni{1} $A$-paths and homotopy
\ssni{1.1} $A$-paths
\ssni{1.2} $A$-paths and connections
\ssni{1.3} Homotopy of $A$-paths
\ssni{1.4} Representations and $A$-paths
\sni{2} The Weinstein groupoid
\ssni{2.1} The groupoid ${\cal G}(A)$
\ssni{2.2} Homomorphisms
\ssni{2.3} The exponential map
\sni{3} Monodromy
\ssni{3.1} Monodromy groups
\ssni{3.2} A second-order monodromy map
\ssni{3.3} Computing the monodromy
\ssni{3.4} Measuring the monodromy
\sni{4} Obstructions to integrability
\ssni{4.1} The main theorem
\ssni{4.2} The Weinstein groupoid as a leaf space
\ssni{4.3} Proof of the main theorem
\sni{5} Examples and applications
\ssni{5.1} Local integrability
\ssni{5.2} Integrability criteria
\ssni{5.3} Tranversally parallelizable foliations
\smallbreak\noindent Appendix A. Flows
\ssni{A.1} Flows and infinitesimal flows
\ssni{A.2} The infinitesimal flow of a section
\smallbreak\noindent References
 
\advance\sectioncount by -1
\section{Introduction} \label{Integrability: Outline} 

This paper is concerned with the general problem of integrability
of geometric structures. The geometric structures we consider are
always associated with local Lie brackets $[~,~]$ on sections of
some vector bundles, or what one calls {\it Lie algebroids}. A Lie
algebroid can be thought of as a generalization of the
tangent bundle, the locus where infinitesimal geometry takes place.
Roughly speaking, the general integrability problem asks for the
existence of a ``space of arrows'' and a product which unravels the
infinitesimal structure. These global objects are usually known as
{\it Lie groupoids} (or {\it differentiable groupoids}) and in this
paper we shall give the precise obstructions to integrate a Lie
algebroid to a Lie groupoid. For an introduction to this problem
and a brief historical account we refer the reader to the recent
monograph \cite{CaWe}. More background material and further references
can be found in \cite{Mack1}, \cite{Mack2}.

To describe our results, let us start by recalling that a {\it Lie
algebroid} over a manifold $M$ consists of a vector bundle $A$ over
$M$, endowed with a Lie bracket $[~,~]$ on the space of sections
$\Gamma(\A)$, together with a bundle map $\#:A\to TM$, called
{\it the anchor}. One requires the induced map $\#:
\Gamma(\A)\to\X^1(M)$ (\footnote{We denote by $\Omega^r(M)$ and
$\X^r(M)$, respectively, the spaces of differential $r$-forms and
$r$-multivector fields on a manifold $M$.  If $E$ is a bundle over
$M$, $\Gamma(E)$ will denote the space of global sections.}) to be a
Lie algebra map, and also the Leibniz identity
$$  [\al,f\beta]=f[\al,\beta]+\#\al(f)\beta,$$ 
to hold, where the vector field $\#\al$ acts on $f$.

For any $x\in M$, there is an induced Lie bracket on
$$  \gg_{x}\equiv \Ker(\#_{x})\subset \A_x   $$ 
which makes it into a Lie algebra. In general, the dimension of
$\gg_x$ varies with~$x$. The image of $\#$ defines a smooth
generalized distribution in $M$, in the sense of Sussmann
(\cite{Suss}), which is integrable.  When we restrict to a leaf $L$
of the associated foliation, the $\gg_x$'s are all isomorphic and
fit into a Lie algebra bundle $\gg_{L}$ over $L$ (see
\cite{Mack1}). In fact, there is an induced Lie algebroid
$$  \A_{L}= \A|_{L} $$ 
which is transitive (i.e.\ the anchor is surjective), and $\gg_{L}$
is the kernel of its anchor map. A general Lie algebroid $\A$ can
be thought of as a singular foliation on $M$, together with
transitive algebroids $\A_{L}$ over the leaves $L$, glued in some
complicated way.
\vskip 10 pt

Before giving the definitions of Lie groupoids and integrability of
Lie algebroids, we illustrate the problem by looking at the following
basic examples:

\begin{itemize}
\item For algebroids over a point (i.e.\ Lie algebras), the
integrability problem is solved by Lie's third theorem on the
integrability of (finite-dimensional) Lie algebras by Lie groups;
\item For algebroids with zero anchor map (i.e.\ bundles of Lie
algebras), it is Douady-Lazard \cite{DoLa} extension of Lie's third
theorem which ensures that the Lie groups integrating each Lie algebra
fiber fit into a smooth bundle of Lie groups;
\item For algebroids with injective anchor map (i.e.\ involutive
distributions $\F\subset TM$), the integrability problem is solved
by Frobenius' integrability theorem.
\end{itemize}
Other fundamental examples come from \'Elie Cartan's infinite continuous
groups (Singer and Sternberg, \cite{Sing}), the integrability of
infinitesimal actions of Lie algebras on manifolds (Palais,
\cite{Palais}), abstract Atiyah sequences (Almeida and Molino,
\cite{Mol}; Mackenzie, \cite{Mack1}), of Poisson manifolds (Weinstein,
\cite{Wein}) and of algebras of vector fields (Nistor,
\cite{Ni}). These, together with various other examples will be
discussed in the forthcoming sections.

Let us look closer at the most trivial example. A vector field
$X\in\X^1(M)$ is the same as a Lie algebroid structure on the trivial
line bundle $\L\to M$: the anchor is just multiplication by
$X$, while the Lie bracket on $\Gamma(\L)\simeq C^{\infty}(M)$ is
given by $[f, g]= X(f)g- fX(g)$. The integrability result here
states that a vector field is integrable to a local flow. It may be
useful to think of the flow $\Phi_{X}^{t}$ as a collection of
arrows $x\rmap\Phi_{X}^{t}(x)$ between the different points of the
manifold, which can be composed by the rule
$\Phi_{X}^{t}\Phi_{X}^{s}=\Phi_{X}^{s+t}$. The points which can be
joined by such an arrow with a given point $x$ form the orbit of
$\Phi_{X}$ (or the integral curve of $X$) through $x$.

The general integrability problem is similar: it asks for the
existence of a ``space of arrows'' and a partially defined
multiplication, which unravels the infinitesimal structure $(A,
[~,~], \#)$. In a more precise fashion, a {\it groupoid} is a
small category $\G$ with all arrows invertible. If the set
of objects (points) is $M$, we say that $\G$ is a groupoid over
$M$. We shall denote by the same letter $\G$ the space of arrows,
and write
$$
\begin{array}{c}
  \G\\
{\scriptstyle\s}\Big\downarrow\Big\downarrow {\scriptstyle \t}\\
  M
\end{array}
$$ 
where $\s$ and $\t$ are the source and target maps. If $g,h\in\G$ the
product $gh$ is defined only for pairs $(g,h)$ in the set of
composable arrows
$$  \G^{(2)}=\set{(g,h)\in\G\times\G| \t(h)=\s(g)},$$ 
and we denote by $g^{-1}\in \G$ the inverse of $g$, and by $1_{x}= x$
the identity arrow at $x\in M$. If $\G$ and $M$ are topological
spaces, all the maps are continuous, and $\s$ and $\t$ are open
surjections, we say that $\G$ is a {\it topological groupoid}. A
{\it Lie groupoid} is a groupoid where the space of arrows $\G$ and
the space of objects $M$ are smooth manifolds, the source and target
maps $\s, \t$ are submersions, and all the other structure maps are
smooth. We require $M$ and the $\s$-fibers $\G(x, -)= \s^{-1}(x)$,
where $x\in M$, to be Hausdorff manifolds, but it is important to
allow the total space $\G$ of arrows to be non-Hausdorff: simple examples
arise even when integrating Lie algebra bundles \cite{DoLa},
while in foliation theory it is well known that the
monodromy groupoid of a foliation is non-Hausdorff if there are
vanishing cycles. For more details see \cite{CaWe}.
\vskip 10 pt

As in the case of Lie groups, any Lie groupoid $\G$ has an associated
Lie algebroid $A= A(\G)$. As a vector bundle, it is the restriction to
$M$ of the bundle $T^{\s}\G$ of $\s$-vertical vector fields on
$M$. Its fiber at $x\in M$ is the tangent space at $1_x$ of the
$\s$-fibers $\G (x, -)= \s^{-1}(x)$, and the anchor map is just the
differential of the target map $\t$. To define the bracket, one shows
that $\Gamma(\A)$ can be identified with $\X_{\rm inv}^{\s}(\G)$, the space
of $\s$-vertical, right-invariant, vector fields on $\G$.  The
standard formula of Lie brackets in terms of flows shows that
$\X_{\rm inv}^{\s}(\G)$ is closed under $[\cdot, \cdot]$. This induces a
Lie bracket on $\Gamma(\A)$, which makes $\A$ into a Lie algebroid.

We say that a Lie algebroid $\A$ is integrable if there exists a Lie
groupoid $\G$ inducing $\A$. The extension of Lie's theory (Lie's
first and second theorem) to Lie algebroids has a promising start.

\nonumproclaim{Theorem {\rm (Lie I)}}
If $\A$ is an integrable Lie algebroid{\rm ,} then there exists a
{\rm (}\/unique\/{\rm )} $\s$\/{\rm -}\/simply connected Lie groupoid integrating $\A$.
\endproclaim

This has been proved in \cite{MoMr} (see also~\cite{Mack1} for the
transitive case).  A different argument, which is just an extension of
the construction of the smooth structure on the universal cover of a
manifold (cf.\ Theorem 1.13.1 in~\cite{DuKo}), will be presented
below. Here $\s$-simply connected means that the $\s$-fibers
$\s^{-1}(x)$ are simply connected. The Lie groupoid in the theorem is
often called the {\it monodromy groupoid} of $\A$, and will be
denoted by $\Mon(\A)$. For the simple examples above, $\Mon(TM)$ is
the homotopy groupoid of $M$, $\Mon(\F)$ is the monodromy groupoid of
the foliation $\F$, while $\Mon(\gg)$ is the unique simply-connected
Lie group integrating $\gg$.

The following result is standard (we refer to ~\cite{MaXu}, \cite{MoMr},
although the reader may come across it in various other places).
See also Section \ref{section:Weinstein groupoid} below.

\nonumproclaim{Theorem {\rm (Lie II)}}
Let $\phi: A\to B$ be a morphism of integrable Lie algebroids{\rm ,} and
let $\G$ and $\H$ be integrations of $\A$ and $B$. If $\G$ is
$\s$\/{\rm -}\/simply connected{\rm ,} then there exists a {\rm (}\/unique\/{\rm )} morphism of Lie
groupoids $\Phi: \G\to \H$ integrating $\phi$.
\endproclaim

In contrast with the case of Lie algebras or foliations, there is no
Lie's third theorem for general Lie algebroids.  Examples of
nonintegrable Lie algebroids are known (we will see several of them in
the forthcoming sections) and, up to now, no good explanation for this
failure was known. For transitive Lie algebroids, there is a cohomological
obstruction due to Mackenzie (\cite{Mack1}), which may be regarded as an
extension
to non-abelian groups of the Chern class of a circle bundle, and which
gives a necessary and sufficient criterion for integrability. Other
various integrability criteria one finds in the literature are
(apparently) nonrelated: some require a nice behavior of the Lie algebras
$\gg_x$, some require a nice topology of the leaves of the induced
foliation, and most of them require regular algebroids. A good
understanding of this failure should shed some light on the following
questions:

\begin{itemize}
\item {\it Is there a {\rm (}\/computable\/{\rm )} obstruction to the integrability
of Lie algebroids}\/?
\item {\it Is the integrability problem a local one}\/?
\item {\it Are Lie algebroids locally integrable}\/?
\end{itemize}

In this paper we provide answers to these questions. We show
that the obstruction to integrability comes from the relation between
the topology of the leaves of the induced foliation and the Lie
algebras defined by the kernel of the anchor map.
\vskip 10 pt

We will now outline our integrability result. Given an algebroid $\A$ and
$x\in M$, we will construct certain (monodromy) subgroups
$N_{x}(\A) \subset \A_x,$
which lie in the center of the Lie algebra $\gg_x= {\rm Ker}(\#_{x})$:
they consist of those elements $v\in Z(\gg_x)$ which are
{\it homotopic} to zero (see \S \ref{section:A-paths}). As we shall
explain, these groups arise as the image of a second-order monodromy map
$$  \partial: \pi_{2}(L_x)\to \G(\gg_x),$$ 
which relates the topology of the leaf $L_x$ through $x$ with the
simply connected Lie group $\G(\gg_x)$ integrating the Lie algebra
$\gg_x= {\rm Ker}(\#_{x})$. {From} a conceptual point of view, the monodromy map
can be viewed as an analogue of a boundary map of the homotopy long exact
sequence of a fibration (namely $0\to \gg_{L_x}\to A_{L_x}\to
TL_x\to 0$). In order to measure the discreteness of the
groups $N_{x}(\A)$ we let
$$  r(x)= d(0, N_{x}(\A)-\{0\}),$$ 
where the distance is computed with respect to a (arbitrary)  norm on
the vector bundle $\A$. Here we adopt the convention
$d(0,\emptyset)=+\infty$.
We will see that $r$ is {\it not} a continuous function. Our main result
is:

\nonumproclaim{Theorem {\rm (Obstructions to Lie III)}} For a Lie algebroid $\A$ over
$M${\rm ,}
the fol\-low\-ing are equivalent\/{\rm :}
\begin{itemize}
\item[{\rm (i)}] $\A$ is integrable{\rm ;}
\item[{\rm (ii)}] For all $x\in M$, $N_{x}(A)\subset \A_x$ is discrete and
$\liminf_{y\to x}r(y)>0${\rm .}
\end{itemize}
\endproclaim

We stress that these obstructions are computable in many examples.
First of all, the definition of the monodromy map is explicit.
Moreover, given a splitting $\sigma:TL\to A_L$ of $\#$ with
$Z(\gg_L)$-valued curvature 2-form $\Omega_{\sigma}$, we will see that
$$  N_{x}(A)= \{ \int_{\gamma} \Omega_{\sigma}:
\gamma\in \pi_{2}(L, x)\}\subset Z(\gg_x).$$ 
With this information at hand the reader can already jump to the examples
(see \S\S 3.3, 3.4, 4.1 and \ref{Examples/Applications}).
\vskip 10 pt

As is often the case, the main theorem is just an instance of a more
fruitful
approach. In fact, we will show that a Lie algebroid $A$ always
admits an ``integrating'' {\it topological} groupoid $\G(\A)$. Although
it is not always smooth (in general it is only a leaf space), it does behave
like a Lie groupoid. This immediately implies the integrability of Lie
algebroids by ``local Lie groupoids'', a result which has been assumed to
hold 
since the original works of Pradines in the 1960's.

The main idea of our approach is as follows: Suppose $\pi:\A\to M$
is a Lie algebroid which can be integrated to a Lie groupoid
$\G$. Denote by $P(\G)$ the space of {\it $\G$-paths}, with the
$C^2$-topology:
$$  P(\G)=\set{g: [0,1]\to \G|~ g\in C^2,\ \s(g(t))=x,\ g(0)= 1_{x}}$$ 
(paths lying in $\s$-fibers of $\G$ starting at the identity).
Also, denote by $\sim$ the equivalence relation defined by
$C^1$-homotopies in $P(\G)$ with fixed end-points. Then we have a
standard description of the monodromy groupoid as
$$  \Mon(\A)= P(\G)/\sim.$$ 
The source and target maps are the obvious ones, and for two paths $g, g'\in
P(\G)$ which are composable (i.e.\ $\t(g(1))=\s(g'(0))$) we define
$$ 
g'\cdot g(t)\equiv \left\{
\begin{array}{ll}
g(2t),\qquad& 0\le t\le \frac{1}{2}\\ \\
g'(2t-1)g(1),\qquad & \frac{1}{2}< t\le 1 .
\end{array}
\right. 
$$ 
Note that any element in $P(\G)$ is equivalent to some $g(t)$ with
derivatives vanishing at the end-points, and if $g$ and $g'$ have this
property, then $g'\cdot g\in P(\G)$. Therefore, this multiplication is
associative up to homotopy, so we get the desired multiplication on
the quotient space which makes $\Mon(\A)$ into a (topological)
groupoid. The construction of the smooth structure on $\Mon(\A)$ is
similar to the construction of the smooth structure on the universal
cover of a manifold (see e.g.\ Theorem 1.13.1 in \cite{DuKo}).

Now, any $\G$-path $g$ defines an $\A$-path $a$, i.e.\ a curve
$a:I\to \A$ defined on the unit interval $I= [0, 1]$, with the
property that
$$  \#a(t)= \frac{d}{dt} \pi(a(t)).$$ 
The $\A$-path $a$ is obtained from $g$ by differentiation and right
translations. This defines a bijection between $P(\G)$ and the set
$P(\A)$ of $\A$-paths and, using this bijection, we can transport
homotopy of $\G$-paths to an equivalence relation ({\it homotopy}) of
$\A$-paths.  Moreover, this equivalence can be expressed using the
infinitesimal data only (\S \ref{section:A-paths}, below).  It
follows that a monodromy type groupoid $\G(\A)$ can be constructed
without any integrability assumption. This construction of $\G(\A)$,
suggested by Alan Weinstein, in general only produces a
{\it topological} groupoid (\S \ref{section:Weinstein
groupoid}).  Our main task will then be to understand when does the
Weinstein groupoid $\G(\A)$ admit the desired smooth structure, and
that is where the obstructions show up. We first describe the second-order monodromy map which encodes these obstructions (\S
\ref{section:monodromy}) and we then show that these are in fact the
only obstructions to integrability (\S 
\ref{section:integrability}). In the final section, we derive the
known integrability criteria from our general result and we give
two applications.

%
\demo{Acknowledgments}
The construction of the groupoid $\G(A)$ was suggested to us by
Alan Weinstein, and is inspired by a ``new'' proof of Lie's third
theorem in the recent monograph \cite{DuKo} by Duistermaat and
Kolk.  We are indebted to him for this suggestion as well as many
comments and discussions. \pagebreak The same type of construction, for the
special case of Poisson manifolds, appears in the work of Cattaneo
and Felder \cite{CaFe}. Though they do not discuss integrability
obstructions, their paper was also a source of inspiration for the
present work.

We would also like to express our gratitude for
additional comments and discussions to Ana Cannas da Silva, Viktor
Ginzburg, Kirill Mackenzie, Ieke Moerdijk, Janez Mr\v{c}un and James
Stasheff.
\enddemo

\section{$\A$-paths and homotopy}    %
\label{section:A-paths}                      %

In this section $\A$ is a Lie algebroid over $M$, $\#:\A\to TM$
denotes the anchor, and $\pi:\A\to M$ denotes the projection.

In order to construct our main object of study, the groupoid
$\G(\A)$ that plays the role of the monodromy groupoid $\Mon(\A)$
for a general (nonintegrable) algebroid, we need the appropriate
notion of paths on $\A$. These are known as {\it $A$-paths} (or
admissible paths) and we shall discuss them in this section.

\demo{{\rm 1.1.} $\A$-paths}
We call a $C^1$ curve $a:I\to \A$ an {\it $\A$\/{\rm -}\/path} if
$$ \#a(t)= \frac{d}{dt} \gamma(t),$$ 
where $\gamma(t)= \pi(a(t))$ is the base path (necessarily of class $C^2$).
We let $P(\A)$ denote the space of $A$-paths, endowed with the
topology of uniform convergence.

We emphasize that this is the right notion of paths in the world of
algebroids. {From} this point of view, one should view $a$ as a bundle map
$$  a\,dt: TI\to A $$ 
which covers the base path $\gamma: I\to M$ and this gives a
algebroid morphism $TI\to A$.

Obviously, the base path of an $\A$-path sits inside a leaf $L$ of
the induced foliation, and so can be viewed as an $\A_{L}$-path. The
key remark is: \enddemo

\proclaim{Proposition}\label{A-G-paths}
If $\G$ integrates the Lie algebroid $A${\rm ,} then there is a
homeomorphism $D^{R}: P(\G)\to P(\A)$ between the space of $\G$\/{\rm -}\/paths{\rm ,}
and the space of $\A$\/{\rm -}\/paths {\rm (}$D^{R}$ is called the differentiation of
$\G$\/{\rm -}\/paths{\rm ,} and its inverse is called the integration of $\A$\/{\rm -}\/paths\/{\rm .)}
\endproclaim 

\demo{Proof} 
Any $\G$-path $g:I\to \G$ defines an $A$-path $D^{R}(g):I\to
A$ by the formula
$$  (D^{R}g)(t)= (dR_{g(t)^{-1}})_{g(t)} \dot{g}(t) \ ,$$ 
where, for $h: x\to y$ an  arrow in $\G$, $R_{h}: \s^{-1}(y)\to s^{-1}(x)$
is the right multiplication by $h$.  Conversely, any $A$-path $a$
arises in this way, by integrating (using \textsc{Lie II}) the Lie
algebroid morphism $TI\to A$ defined by $a$. Finally, notice that
any Lie groupoid homomorphism $\phi:I\times I\to \G$ from the pair
groupoid into $\G$, is of the form $\phi(s,t)=g(s)g^{-1}(t)$ for some
$\G$-path $g$.

A more explicit argument, avoiding \textsc{Lie II}, and which also
shows that the inverse of $D^{R}$ is continuous, is as follows.
Given $a$, we choose a time-dependent section $\alpha$ of $\A$
extending $a$, i.e.\  so that
$$  a(t)= \alpha(t, \gamma(t)).$$ 
If we let $\varphi_{\alpha}^{t, 0}$ be the flow of the right-invariant
vector field that corresponds to $\alpha$, then $g(t)=
\varphi_{\alpha}^{t, 0}(\gamma(0))$ is the desired $\G$-path.  Indeed,
right-invariance guarantees that this flow is defined for all
$t\in[0,1]$ and also implies that
\vglue12pt
\hfill $ (D^{R}g)(t)=(dR_{g(t)^{-1}})_{g(t)}(\alpha(t, g(t)))=
                                \alpha(t, \gamma(t))= a(t).$ 
\enddemo

\demo{{\rm 1.2.} $\A$-paths and connections} Given an $\A$-connection on a vector bundle $E$ over $M$, most of the
classical constructions (which we recover when $\A= TM$) extend to Lie
algebroids, provided we use $\A$-paths. This is explained in detail in
\cite{Fer1}, \cite{Fer2}, and here we recall only the results we need.

An {\it $\A$\/{\rm -}\/connection} on a vector bundle $E$ over $M$
can be defined by an\break $A$-derivative operator
$\Gamma(A)\times \Gamma(E)\to \Gamma(E)$,
$(\alpha, u)\mapsto \nabla_{\alpha}u$
satisfying $\nabla_{f\alpha}u= f\nabla_{\alpha}u$, and
$\nabla_{\alpha}(fu)= f\nabla_{\alpha}u+ \#\alpha(f) u$.
The curvature of $\nabla$ is given by the usual formula
$$  R_{\nabla}(\alpha, \beta)= [\nabla_{\alpha}, \nabla_{\beta}]-
             \nabla_{[\alpha, \beta]},$$ 
and $\nabla$ is called flat if $R_{\nabla}=0$. For an $A$-connection
$\nabla$ on the vector bundle~$A$, the torsion of $\nabla$ is also defined
as usual by:
$$  T_{\nabla}(\alpha,\beta)=\nabla_{\alpha}\beta-
              \nabla_{\beta}\alpha-[\alpha,\beta].$$ 

Given an $\A$-path $a$ with base path $\gamma: I\to M$,
and $u:I\to E$ a path in $E$ above $\gamma$, then
the derivative of $u$ along $a$, denoted $\nabla_{a}u$,
is defined as usual: choose a time-dependent section $\xi$ of $E$
such that $\xi(t, \gamma(t))= u(t)$, then
$$  \nabla_{a}u(t)= \nabla_{a}\xi^{t}(x)+ \frac{d\xi^{t}}{dt}(x),
\ \ {\rm at}\ \ x= \gamma(t)\ .$$ 
One has then the notion of parallel transport along $a$, denoted
$T_{a}^{t}: E_{\gamma(0)}\to E_{\gamma(t)}$, and for the special
case $E= \A$, we can talk about the geodesics of $\nabla$.
Geodesics are $\A$-paths $a$ with the property that
$\nabla_{a}a(t)= 0$. Exactly as in the classical case, one has existence
and uniqueness of geodesics with given initial base point $x\in M$
and  \pagebreak ``initial speed'' $a_{0}\in A_{x_0}$. \enddemo

\numbereddemo{{E}xample}\label{holonomy}
If $L$ is a leaf of the foliation induced by $\A$, then $\gg_{L}=
{\rm Ker}(\# |_{L})$ carries a flat $A_L$-connection defined by
$\nabla_{\alpha}\beta= [\alpha, \beta]$. In particular, for any
$\A$-path $a$, the induced parallel transport defines a linear map,
called {\it the linear holonomy} of $a$,
$$  \Hol(a): \gg_{x} \to \gg_{y} \ ,$$ 
where $x, y$ are the initial and the end-point of the base path. For
more on linear holonomy we refer to \cite{Fer1}.
\enddemo

Most of the connections that we will use are induced by a standard
$TM$-connection $\nabla$ on the vector bundle $\A$. Associated with
$\nabla$ there is an obvious $\A$-connection on the vector bundle $\A$
$$  \nabla_{\alpha}\beta\equiv \nabla_{\#\alpha}\beta.$$ 
A bit more subtle are the following two $\A$-connections on $\A$ and
on $TM$, respectively  (see \cite{Cra2}):
$$  \overline{\nabla}_{\alpha}\beta\equiv
                \nabla_{\#\beta}\alpha+ [\alpha, \beta], \qquad
   \overline{\nabla}_{\alpha}X\equiv \# \nabla_{X}\alpha+
   [\#\alpha, X] .$$ 
Note that
$\overline{\nabla}_{\alpha}\#\beta=\#\overline{\nabla}_{\alpha}\beta$,
so in the terminology of \cite{Fer1} this means that
$\overline{\nabla}$ is a basic connection on $\A$. These connections
play a fundamental role in the theory of characteristic classes (see
\cite{Cra}, \cite{Cra2}, \cite{Fer1}).

\demo{{\rm 1.3.} Homotopy of $\A$-paths}  As we saw above, if $A$ is integrable, $A$-paths are in a bijective
correspondence with $\G$-paths.  Let us see now how one can transport
the notion of homotopy to $P(\A)$, so that it only uses the
infinitesimal data (i.e., Lie algebroid data).

Let us fix
$$ a_{\epsilon}(t)= a(\epsilon, t): I\times I\to A$$ 
a {\it variation of $\A$-paths}, that is a family of $\A$-paths
$a_{\epsilon}$ which is of class $C^2$ on $\epsilon$, with the
property that the base paths $\gamma_{\epsilon}(t)= \gamma
(\epsilon, t): I\times I\to M$ have fixed end-points.  If $\A$ came
from a Lie groupoid $\G$, and $a_{\epsilon}$ came from $\G$-paths
$g_{\epsilon}$, then $g_{\epsilon}$ would not necessarily give a homotopy
between $g_{0}$ and $g_{1}$, because the end-points
$g_{\epsilon}(1)$ may vary. The following lemma describes two
distinct ways of controlling the variation $\frac{d}{d\epsilon}
g_{\epsilon}(1)$: one way uses a connection on $A$, and the other
uses flows of sections of a $A$ (see Appendix A). They both depend
only on infinitesimal data.
\enddemo

\proclaim{Proposition}
\label{equivalence}
Let $\A$ be an algebroid and $a= a_{\epsilon}$ a variation of $\A$\/{\rm -}\/paths.
\begin{itemize}
\item[{\rm (i)}] If $\nabla$ is a  $TM$-connection on $\A$ with torsion
$T_{\nabla}${\rm ,} then the solution $b=b(\epsilon, t)$ of the differential
equation \end{itemize}
\begin{equation}\label{diffeq}
\partial_{t}b- \partial_{\epsilon} a= T_{\nabla}(a, b),\qquad
b(\epsilon, 0)=0,
\end{equation}
\begin{itemize} \item[]
does not depend on $\nabla$. Moreover{\rm ,} $\#b= \frac{d}{d\epsilon}\gamma$.
\item[{\rm (ii)}] If $\xi_{\epsilon}$ are time-depending sections of $\A$ such
that $\xi_{\epsilon}(t, \gamma_{\epsilon}(t))= a_{\epsilon}(t)${\rm ,}
then $b(\eps,t)$ is given by \end{itemize}
\begin{equation}
\label{variation}
b(\epsilon, t) = \int_{0}^{t}\phi_{\xi_{\epsilon}}^{t, s}
\frac{d \xi_{\epsilon}}{d\epsilon} (s, \gamma_{\epsilon}(s)) ds,
\end{equation}
 \begin{itemize} 
\item[] where $\phi^{t,s}_{\xi_{\epsilon}}$ denotes the flow of the time-dependent section $\xi_{\epsilon}$.
\item[{\rm (iii)}] If $\G$ integrates $\A$ and $g_{\epsilon}$ are the
$\G$\/{\rm -}\/paths satisfying $D^{R}(g_{\epsilon})= a_{\epsilon}${\rm ,} then $b=
D^{R}(g^{t})${\rm ,} where $g^t$ are the paths in $\G$: $\eps\to
g^{t}(\epsilon)=g(\epsilon, t)$.
\end{itemize}
\endproclaim  

This motivates the following definition:

\numbereddemo{{D}efinition}
We say that two $\A$-paths $a_0$ and $a_1$ are equivalent
(or homotopic), and write $a_0\sim a_1$, if there exists a variation
$a_{\epsilon}$ with the property that $b$ insured by Proposition
 \ref{equivalence} satisfies $b(\epsilon, 1)= 0$ for all $\epsilon \in I$.
\enddemo

When $\A$ admits an integration $\G$, then the isomorphism $D^{R}:
P(\G)\to P(\A)$ of Proposition \ref{A-G-paths} transforms the usual
homotopy into the homotopy of $\A$-paths.  Note also that, as
$\A$-paths should be viewed as algebroid morphisms, the pair $(a, b)$
defining the equivalence of $\A$-paths should be viewed as a true
homotopy
$$  a dt + b d\epsilon  : TI\times TI \to A $$ 
in the world of algebroids. In fact, equation (\ref{diffeq}) is just
an explicit way of saying that this is a morphism of Lie algebroids
(see \cite{HiMa}).

\demo{Proof  of Proposition {\rm \ref{equivalence}}} Obviously, (i)
follows from 
(ii). To prove (ii), let $\xi_{\epsilon}$ be as in the statement, and
let $\eta$ be given by
$$  \eta (\epsilon, t, x)= \int_{0}^{t}
\phi_{\xi_{\epsilon}}^{t, s} \frac{d\xi_{\epsilon}}{d\epsilon}(s,
\Phi_{\#\xi_{\epsilon}}^{s, t}(x)) ds \in A_x .$$ 
We may assume that $\xi_{\epsilon}$ as compact support.  We note that $\eta$ coincides with the solution of the
equation
\begin{equation}
\label{eq.equiv}
\frac{d\eta}{dt}- \frac{d\xi}{d\epsilon}= [\eta, \xi]\ ,
\end{equation}
with $\eta(\epsilon, 0)= 0$. Indeed, since
$$  \eta (\epsilon, t, -)= \int_{0}^{t}
(\phi_{\xi_{\epsilon}}^{s, t})^{*}(\frac{d\xi_{\epsilon}^{s}}{d\epsilon}) ds
\in \Gamma(A),$$ 
equation (\ref{eq.equiv}) immediately follows from the basic
formula (\ref{Lie-flows}) for flows. Also, $X=\#\xi$ and $Y=\#\eta$
satisfy a similar equation on $M$, and since we have
$X(\epsilon,t,\gamma_{\epsilon}(t))= \frac{d\gamma}{dt}$, it
follows that $Y(\epsilon, t, \gamma_{\epsilon}(t))=
\frac{d\gamma}{d\epsilon}$. In other words, $b(\epsilon, t)=
\eta(\epsilon, t, \gamma(\epsilon, t))$ satisfies $\#b=
\frac{d\gamma}{d\epsilon}$. We now have
$$  \partial_{t}b= \nabla_{\frac{d\gamma}{dt}}\eta+ \frac{d\eta}{dt}=
\nabla_{\#\xi}\eta+ \frac{d\eta}{dt} $$ 
at $x= \gamma_{\epsilon}(t)$. Subtracting from this the similar formula
for $\partial_{\epsilon}a$ and using (\ref{eq.equiv}) we get
$$  \partial_{t}b- \partial_{\epsilon}a= \nabla_{\#\xi}\eta-
\nabla_{\#\eta}\xi+ [\eta, \xi]= T_{\nabla}(\xi,\eta) .$$ 

We are now left proving (iii). Assume that $\G$ integrates $\A$ and
$g_{\epsilon}$ are the $\G$-paths satisfying $D^{R}(g_{\epsilon})=
a_{\epsilon}$. The formula of variation of parameters applied to the
right-invariant vector field $\xi_{\epsilon}$ shows that
\begin{eqnarray*}
\frac{\partial g(\eps,t)}{\partial \eps}&=&
\int_0^t (d\varphi^{t,s}_{\xi_{\epsilon}})_{g(\eps,s)}
\frac{d\xi_{\epsilon}^s}{d\eps}(g(\eps,s))ds\\ &=&
(dR_{g(\eps,t)})_{\gamma_\eps(t)}
\int_0^t\phi^{t,s}_{\xi_{\epsilon}}
\frac{d\xi_{\epsilon}^s}{d\eps}(\gamma_\eps(s))ds.
\end{eqnarray*}
But then:
\vglue12pt
\hfill ${\displaystyle
D^R(g^t)=\int_0^t\phi^{t,s}_{\xi_{\epsilon}}
\frac{d\xi_{\epsilon}^s}{d\eps}(\gamma_\eps(s))ds=b(\eps,t).
}$ 
\enddemo

The next lemma gives elementary properties of homotopies of
$A$-paths:

\proclaim{Lemma}
\label{homotopy properties}
Let $\A$ be a Lie algebroid.
\begin{itemize}
\item[{\rm (i)}] If $\tau: I\to I${\rm ,} with $\tau(0)=
0${\rm ,} $\tau(1)= 1$ is a smooth change of parameter{\rm ,} then any
$\A$\/{\rm -}\/path $a$ is equivalent to its reparametrization
$a^{\tau}(t)\equiv\tau'(t)a(\tau(t))$.
\item[{\rm (ii)}] Any $\A$\/{\rm -}\/path $a_{0}$ is equivalent to a smooth {\rm (}\/i.e.\ of
class $C^{\infty}${\rm )} $\A$\/{\rm -}\/path.
\item[{\rm (iii)}] If two smooth $\A$\/{\rm -}\/paths $a_{0}, a_{1}$ are equivalent{\rm ,} then
there
exists a smooth homotopy between them.
\end{itemize}

\endproclaim 

\demo{Proof} 
To prove (i), we consider the variation
$$  a_{\epsilon}(t)=
((1-\epsilon)+\epsilon\tau'(t))a((1-\epsilon)t+\epsilon\tau(t)) $$ 
and we check that the associated $b$ satisfies $b(\epsilon, 1)= 0$.
In fact, one can compute by any of the methods of
Proposition \ref{equivalence}:
$$  b(\epsilon,t)= (\tau(t)- t)a((1-\epsilon)t+\epsilon\tau(t)).$$ 
For example, if we let $\al$ be a time-dependent section which extends
the path $a$, and define a 1-parameter family of time-dependent
sections $\xi_{\epsilon}$ by:
$$  \xi_{\epsilon}(t, x)=
((1-\epsilon)+\epsilon\tau'(t))\al((1-\epsilon)t+\epsilon\tau(t),x),$$ 
then $\xi_{\epsilon}$ extends $a_{\epsilon}$ and the family
$$  \eta (\epsilon, t, x)= (\tau(t)- t) \al((1-\epsilon)t+\epsilon\tau(t), x)
$$ 
satisfies (\ref{eq.equiv}). Hence, we must have $b(\epsilon,
t)=\eta(\epsilon, t, \gamma(\epsilon, t))$ as claimed.

For (ii), note that from the similar claim for ordinary paths on
manifolds (see e.g.\ Theorem $1.13.1$ in \cite{DuKo}), we can find
a $C^{r}$-homotopy $\gamma_{\epsilon}$ between the base path
$\gamma_{0}$ of $a_0$ and a smooth path $\gamma_{1}$. Also, we
can do it so that $\gamma_{\epsilon}$ stays in the same leaf $L$ as
$\gamma_{0}$, and so that $\gamma_{\epsilon}(t)$ is smooth in the
domain $t\in [0, 1]$, $\epsilon\in [c, 1]$ for some constant $0<c<
1$. We now choose a smooth splitting $\sigma: TL\to A|_{L}$ of the
anchor map, and put $b(\epsilon , t)=
\sigma(\frac{d}{d\epsilon}\gamma_{\epsilon}(t))$. Let $a$ be the
solution of the differential equation (\ref{diffeq}), with the
initial conditions $a(0, t)= a_{0}(t)$. Clearly $a$ is smooth on
the domain on which $b$ is; hence it defines a homotopy between
$a_{0}$ and the smooth $\A$-path $a_{1}$. Part (iii) is just a
degree-one higher version of part (ii), and can be proved similarly,
replacing the path $a_0$ by the given homotopy between $a_0$ and
$a_1$ (a similar argument will be presented in detail in the proof
of Proposition \ref{monodromy-seq}).
\enddemo

1.4. {\it Representations and $A$\/{\rm -}\/paths.} A flat $A$-connection on a vector bundle $E$ defines a
{\it representation} of $A$ on $E$. The terminology is inspired by
the case of Lie algebras. There is also an obvious notion of
representation of a Lie groupoid $\G$: this is a vector bundle $E$
over the space $M$ of objects, together with smooth linear actions $g:
E_{x}\to E_{y}$ defined for arrows $g$ from $x$ to $y$ in $\G$,
satisfying the usual identities. By differentiation, any such
representation becomes a representation of the Lie algebroid $\A$ of
$\G$ (see e.g.\ \cite{Cra}, \cite{HiMa}). Moreover, when $\G= \Mon(\A)$ is
the unique $\s$-simply connected Lie groupoid integrating $\A$, this
construction induces a bijection
$$  \Rep(\Mon(\A)) \cong \Rep(\A) $$ 
between the (semi-rings of equivalence classes of) representations.
This is explained in \cite{Cra}, \cite{Ginz}, using the integrability of actions of
\cite{MoMr},
but it follows also from our construction of $\G(\A)$ (see next section)
since we have:

\proclaim{Proposition}
\label{hol-hom}
If $a_0$ and $a_1$ are equivalent $\A$\/{\rm -}\/paths from $x$ to $y$. Then
for any representation $E$ of $\A${\rm ,} parallel transports $E_x\to
E_y$ along $a_0$ and $a_1$ coincide.
\endproclaim

\demo{Proof}  We first claim that for any $\A$-connection $\nabla$ on $E$,
 and homotopy $adt+ bd\epsilon$  between $a_0$ and $a_1$, we have:
$$ 
\nabla_{a_{\epsilon}}\nabla_{b_{t}}u- \nabla_{b_{t}}\nabla_{a_{\epsilon}}u=
R_{\nabla}(a, b) u
$$ 
for all paths $u: I\times I\to E$ above $\gamma(\epsilon, t)$. To
see this, let us assume that $\xi$, $\eta$ are as in the proof
of Proposition \ref{equivalence}, and let $s$ be a family of
time-dependent sections of $E$ so that $u(\epsilon, t)=
s(\epsilon, t, \gamma(\epsilon, t))$. Then
$$  \nabla_{b_{t}}u= \nabla_{\eta}s+ \frac{d s}{d\epsilon} $$ 
at $x= \gamma(\epsilon, t)$.  Hence
$$ 
\nabla_{a_{\epsilon}}\nabla_{b_{t}}u
=\nabla_{\xi}\nabla_{\eta}s+ \nabla_{\xi}(\frac{d s}{d\epsilon})+
       \nabla_{\eta}(\frac{d s}{dt})+ \frac{d^2 s}{d\epsilon dt}
       +\nabla_{\frac{d\eta}{dt}}s.
$$ 
Subtracting the analogous formula for
$\nabla_{b_{t}}\nabla_{a_{\epsilon}}u$ and using (\ref{eq.equiv}), we have
proved the claim.

When $\nabla$ is flat, this formula applied to $u(\epsilon, t)=
T_{a_{\epsilon}}^{t}(u_0)$, where $T_{a_{\epsilon}}^{t}$ denotes
parallel transport, gives $\nabla_{a_{\epsilon}}\nabla_{b_{t}}u=
0$. But $\nabla_{b_{t}}u= 0$ at $t= 0$, hence $\nabla_{b_{t}}u= 0$
for all $t$'s. Since $u(0, t)= T_{a_0}^{t}(u_0)$ it follows that
$u(\epsilon, t)= T_{b_t}^{\epsilon} T_{a_0}^{t}(u_0)$. Therefore
$T_{a_{\epsilon}}^{t}= T_{b_{t}}^{\epsilon}T_{a_0}^{t}$,
for all $\epsilon, t$ and, in particular, for $\epsilon= t= 1$ we get
$T_{a_1}^{1}= T_{a_0}^{1}$.
\enddemo

Recalling the notion of linear holonomy (cf.~Example
\ref{holonomy}) we have:

\proclaim{{C}orollary}
\label{cor:linear:holonomy}
If $a_0$ and $a_1$ are equivalent $\A$\/{\rm -}\/paths from $x$ to $y${\rm ,} they
induce the same linear holonomy maps
$$  \Hol(a_0)=\Hol(a_1): \gg_x\to \gg_y.$$ 
\endproclaim 
\vglue-30pt
 
\section{The Weinstein groupoid}         %
\label{section:Weinstein groupoid}       %
 \vglue-6pt

We are now ready to define the Weinstein groupoid $\G(\A)$ of a
general Lie algebroid, which in the integrable case will be the
unique $\s$-simply connected groupoid integrating $\A$.

\demo{{\rm 2.1.} The groupoid $\G(\A)$}  Let $a_0$, $a_1$ be two composable $\A$-paths, i.e.\ so that
$\pi(a_{0}(1))= \pi(a_{1}(0))$. We define their concatenation
$$ 
a_{1}\odot a_{0}(t)\equiv \left\{
\begin{array}{ll}
2a_{0}(2t),\qquad& 0\le t\le \frac{1}{2}\\ \\
2a_{1}(2t-1),\qquad & \frac{1}{2}< t\le 1.
\end{array}
\right.
$$ 
This is essentially the multiplication that we need. However,
$a_{1}\odot a_{0}$ is only {\it piecewise} smooth. One way around this
difficulty is allowing for $\A$-paths which are {\it piecewise}
smooth. Instead, let us fix a cutoff function $\tau\in
C^\infty(\Rr)$ with the following properties:
\begin{itemize}
\item[(a)] $\tau(t)=1$ for $t\ge 1$ and $\tau(t)=0$ for $t\le 0$;
\item[(b)] $\tau'(t)>0$ for $t\in\, ]0,1[$.
\end{itemize}
\pagebreak \noindent For an $\A$-path $a$ we denote, as above, by $a^{\tau}$ its
reparametrization $a^{\tau}(t)= \tau'(t)a(\tau(t))$. We now define
the multiplication by
$$   a_{1}a_{0}\equiv a_{1}^{\tau}\odot a_{0}^{\tau}\in P(\A).$$ 
According to Lemma \ref{homotopy properties} (i), $a_{0}a_{1}$ is
equivalent to $a_{0}\odot a_{1}$ whenever $a_{0}(1)= a_{1}(0)$.
We also consider the natural structure maps: source and target
$\s,\t: P(\A)\to M$ which map $a$ to $\pi(a(0))$ and $\pi(a(1))$,
respectively, the identity section $\varepsilon:M\to P(\A)$
mapping $x$ to the constant path above $x$, and the inverse $\iota:
P(\A)\to P(\A)$ mapping $a$ to $\overline{a}$ given by
$\overline{a}(t)=-a(1-t)$.
\enddemo
 
\vglue-6pt
\proclaim{Theorem}
\label{Weinstein}
Let $\A$ be a Lie algebroid over $M$. Then the quotient
$$  \G(\A)\equiv P(\A)/\sim $$ 
is a $\s$\/{\rm -}\/simply connected topological groupoid independent of
the choice of cutoff function. Moreover{\rm ,} whenever $\A$ is
integrable{\rm ,} $\G(\A)$ admits a smooth structure which makes it into the
unique $\s$\/{\rm -}\/simply connected Lie groupoid integrating $\A$.
\endproclaim 
\vglue-6pt

{\it Proof}.
If we take the maps on the quotient induced from the structure maps
defined above, then $\G(\A)$ is clearly a groupoid. Note that the
multiplication on $P(\A)$ was defined so that, whenever $\G$
integrates $\A$, the map $D^{R}$ of Proposition \ref{A-G-paths}
preserves multiplications. Hence the only thing we still have to
prove is that $\s,\t: \G(\A)\to M$ are open maps.

Given $D\subset \G(\A)$ open, we will show that its saturation
$\tilde{D}$ with respect to the equivalence relation $\sim$ is still open.
This follows from the fact, to be shown later in Theorem
\ref{G-as-fol}, that the equivalence relation can be defined by a
foliation on $P(\A)$.

A more direct argument is to show that for any two equivalent
$\A$-paths $a_{0}$ and $a_{1}$, there exists a homeomorphism of $T:
P(\A)\to P(\A)$ such that $T(a)\sim a$ for all $a$'s, and $T(a_{0})=
a_{1}$. To construct such a $T$ we let $\eta=\eta(\epsilon, t)$ be a
family of time-dependent sections of $\A$ which determines the
equivalence $a_{0}\sim a_{1}$ (see Proposition \ref{equivalence}), so
that $\eta(\epsilon, 0)=\eta(\epsilon, 1)= 0$ (we may assume $\eta$
has compact support, so that all the flows involved are everywhere
defined). Given an $\A$-path $b_0$, we consider a time-dependent
section $\xi_0$ so that $\xi_0(t,\gamma_0(t))= b(t)$ and denote by
$\xi$ the solution of equation (\ref{eq.equiv}) with initial condition
$\xi_{0}$.  If we set
$\gamma_{\epsilon}(t)=\Phi_{\#\eta_{t}}^{\epsilon, 0}\gamma_0(t))$ and
$b_{\epsilon}(t)= \xi_{\epsilon}(t, \gamma_{\epsilon}(t))$, then
$T_{\eta}(b_0)\equiv b_1$ is homotopic to $b_0$ via $b_{\epsilon}$,
and maps $a_{0}$ into $a_{1}$.
\hfill\qed\vglue6pt

2.2. {\it Homomorphisms.} Note that, although $\G(\A)$ is not always smooth, in many aspects
it behaves like in the smooth (i.e.\ integrable) case. For
instance, we can call a representation of $\G(\A)$ {\it smooth} if
the action becomes smooth when pull backed to $P(\A)$. Similarly
one can talk about smooth functions on $\G(\A)$, about its tangent
space, etc. This subsection and the next are variations on this
theme.

\proclaim{Proposition}
\label{prop:wenstein:homomorphism} Let $A$ and $B$ be Lie algebroids. Then\/{\rm :}
\begin{itemize}
\item[{\rm (i)}] Every algebroid homomorphism $\phi:A\to B$ determines
a smooth groupoid homomorphism $\Phi:\G(A)\to \G(B)$ of the
associated Weinstein groupoids. If  $A$ and $B$ are integrable{\rm ,}
then $\Phi_*=\phi${\rm ;}
\item[{\rm (ii)}] Every representation $E\in {\rm Rep}(\A)$ determines a smooth
representation of $\G(\A)${\rm ,} which in the integrable case is the
induced representation.
\end{itemize}

\endproclaim 

\demo{Proof} 
For (i) we define $\Phi$ in the only possible way: If $a\in P(A)$
is an $A$-path then $\phi\circ a$ is an $A$-path in $P(B)$.
Moreover, it is easy to see that  if $a_1\sim a_2$ are equivalent
$A$-paths then $\phi\circ a_1\sim\phi\circ a_2$, so we get a
well-defined smooth map $\Phi:\G(A_1)\to\G(A_2)$ by setting
$$  \Phi([a])\equiv[\phi\circ a].$$ 
This map is clearly a groupoid homomorphism.

Part (ii) follows easily from Proposition \ref{hol-hom}.
\enddemo

In particular we see that, as in the smooth case, there is a
bijection between the representations of $\A$ and the (smooth)
representations of $\G(\A)$:
$$  \Rep(\G(\A)) \cong \Rep(\A). $$

\demo{{\rm 2.3.} The exponential map} Assume first that $\G$ is a Lie groupoid integrating $\A$, and
$\nabla$ is a $TM$-connection on $A$. Then the pull-back of
$\nabla$ along the target map $\t$ defines a family of (right-invariant) connections $\nabla_{x}$ on the manifolds
$\s^{-1}(x)$. The associated exponential maps $\Exp_{\nabla_{x}}:
\A_{x}=T_x^{\s}\G\to \s^{-1}(x)$ fit together into a global
exponential map \cite{NiXuWe}
$$  \overline{\Exp}_{\nabla}: A\to \G $$ 
(defined only on an open neighborhood of the zero section). By
standard arguments, $\overline{\Exp}_{\nabla}$ is a diffeomorphism
on a small enough neighborhood of $M$.

Now if $\A$ is not integrable, we still have the exponential map
associated to a connection $\nabla$ on $\A$. It is defined as usual,
so $\Exp_{\nabla}(a)$ is the value at time $t=1$ of the geodesic
($\A$-path) with the initial condition $a$. By a slight abuse of
notation we view it as a map
$$  \Exp_{\nabla}: \A\to P(\A) .$$ 
Of course, $\Exp_{\nabla}$ is only defined on an open neighborhood of
$M$ inside $\A$ consisting of elements whose geodesics are defined for
all $t\in [0,1]$. Passing to the quotient, we have an induced
exponential map
$$  \overline{\Exp}_{\nabla}: A\to \G(\A) .$$ 
For integrable $\A$, this coincides with the exponential map above.

Note that the exponential map we have discussed so far depends on
the choice of the connection. To get an exponential, independent of the
connection, recall (\cite{Mack1}) that an {\it admissible
section} of a Lie groupoid $\G$ is a differentiable map
$\sigma:M\to\G$, such that $\s\circ\sigma(x)=x$ and
$\t\circ\sigma:M\to M$ is a diffeomorphism. Also, each admissible section
$\sigma\in \Gamma(\G)$ determines diffeomorphisms
\begin{eqnarray*}
\G\ni g\mapsto \sigma g\equiv& \sigma(x)g,\qquad \hbox{ where } x=\t(g),\\
\G\ni g\mapsto g \sigma\equiv& g\sigma(y),\qquad \hbox{ where }
\t\circ\sigma(y)=\s(g).
\end{eqnarray*}
Now, each section $\al\in\Gamma(A)$ can be identified with a
right-invariant vector field on $\G$, and we denote its flow by
$\varphi_\al^t$. We define an admissible section $\exp(\al)$ of $\G$
by setting:
$$ \exp(\al)(x)\equiv\varphi_\al^1(x).$$ 
This gives an exponential map $\exp:\Gamma(A)\to \Gamma(\G)$ which, in
general, is defined only for sections $\al$ sufficiently close to the
zero section (e.g., sections with compact support). For more details
see also \cite{Mack1}, \cite{Ni}.

In the nonintegrable case, we can also define an exponential map
$\exp:\Gamma(A)\to \Gamma(\G(A))$ to the admissible smooth sections of
the Weinstein groupoid as follows. First of all notice that
$$  a_{\al}(x)(t)= \al(t, \phi_{\#\al}^{t, 0}(x)) $$ 
defines an $\A$ path $a_{\al}(x)$ for any $x\in M$ and for any time-dependent section $\al$ of $\A$ with flow defined
up to $t=1$ (e.g., if $\al$ is sufficiently close to zero, or if it is compactly
supported).  This defines a smooth map $a_{\al}: M\to P(\A)$. For
$\al\in \Gamma(A)$ close enough to the zero section we set
$$ \exp(\al)(x)=[a_{\al}(x)] .$$ 
Notice that $a=a_{\al}(x)$ is the unique $\A$-path with $a(0)=\al (x)$
and $a(t)=\al (\pi(a(t)))$, for all $t\in I$.

In the integrable case these two constructions coincide. Moreover, for a
general Lie algebroid, we have the following

\proclaim{Proposition}
\label{prop:exponential}
Let $A$ be a Lie algebroid and $\al,\beta\in \Gamma(A)$. Then{\rm ,} as
admissible sections{\rm ,}
$$  \exp(t\al)\exp(\beta)\exp(-t\al)=
\exp(\phi^t_\al \beta)
,$$ 
where $\phi^t_\al$ denotes the infinitesimal flow of $\al$ {\rm (}\/see
Appendix {\rm A).}
\endproclaim 

\demo{Proof} 
First we make the following remark concerning functoriality of $\exp$:
Let $\phi:A_1\to A_2$ be an  isomorphism of Lie algebroids and let
$\Phi:\G(A_1)\to \G(A_2)$ be the corresponding isomorphism of groupoids
(Proposition \ref{prop:wenstein:homomorphism} (i)). If one denotes by
$\tilde{\phi}$ (resp.~$\tilde{\Phi}$) the corresponding homomorphism of
sections (resp.~admissible sections), then we obtain the following
commutative diagram:
$$ 
\begin{array}{ccc}
\Gamma(\G(A_1))&\hskip-5pt\stackrel{\tilde{\Phi}}{\lrar}\hskip-5pt&\Gamma(\G(A_2))\\
{\scriptstyle {\rm exp}} \Big\uparrow&&\Big\uparrow {\scriptstyle {\rm exp}}
\\
\Gamma(A_1)&\hskip-5pt\lower9pt\hbox{$\stackrel{\textstyle\lrar}{\scriptstyle\tilde{\phi}}$}\hskip-5pt&\Gamma(A_2).\end{array}$$

To prove the proposition, it is therefore enough to proof that for
the homomorphism $\Phi^t_\al:\G(A)\to \G(A)$ associated to
$\phi^t_\al:A\to A$ we have:
$$  \Phi^t_\al(g)=\exp(t\al)g\exp(-t\al).$$ 
Equivalently, 
$$  [\phi^t_\al\circ a]=\exp(t\al)[a]\exp(-t\al)$$ 
for any $A$-path $a\in\G(A)$. Now, to prove this, one considers the
variation of
$A$-paths $a_\eps=\exp(-\eps t\al)\cdot (\phi^{\eps t}_\al \circ a)\cdot
\exp(-\eps t\al)$, and checks that this realizes an equivalence of $A$
-paths
using Proposition \ref{equivalence}.
\enddemo

\numbereddemo{{R}emark}
Hence $\G(\A)$ behaves in many respects like a smooth manifold,
even if $\A$ is not integrable. This might be important in various
aspects of noncommutative geometry and its applications to
singular foliations and analysis:
one might expect that the algebras of pseudodifferential operators
and the $C^{*}$-algebra of $\G(\A)$ (see \cite{NiXuWe})
can be constructed even in the nonintegrable case. A related
question is when $\G(\A)$ is a measurable groupoid.
\enddemo

Although the exponential map does exist even in the nonintegrable
case, its injectivity on a neighborhood of $M$ only holds if $\A$ is
integrable. One could say that this is {\it the} difference between
the integrable and the nonintegrable cases, as we will see in the
next sections. However, our main job is to relate the kernel of the
exponential and the geometry of $\A$, and this is the origin of our
obstructions: the monodromy groups described in the next section
consist of the simplest elements which belong to this kernel. It
turns out that these elements are enough to control the entire kernel.

%
\section{Monodromy}                     %
\label{section:monodromy}               %
%

Let $\A$ be a Lie algebroid over $M$, $x\in M$. In this section we
give several descriptions of the (second-order) monodromy
groups of $\A$ at $x$, which control the integrability of $\A$.

\vglue12pt 3.1. {\it Monodromy groups}.
There are several possible ways to introduce the monodromy
groups. Our first description is as follows:

\numbereddemo{{D}efinition}\label{N-groups} We define $N_{x}(A)\subset A_x$ as the subset of the center of ${\frak
g}_x$ formed by those elements $v\in Z(\gg_x)$ with the property that the constant
$\A$-path $v$ is equivalent to the trivial $\A$-path.
\enddemo

Let us denote by $\G(\gg_x)$ the simply-connected Lie group
integrating $\gg_{x}$ (equivalently, the Weinstein groupoid
associated to $\gg_x$). Also, let $\G(A)_x$ be the isotropy
groups of the Weinstein groupoid $\G(\A)$:
$$  \G(A)_x\equiv \s^{-1}(x)\cap \t^{-1}(x) \subset \G(\A)\ .$$ 
Closely related to the groups $N_x(A)$ are the following:

\numbereddemo{{D}efinition}\label{tilde-N-groups} We define $\tilde{N}_x(A)$
as the subgroup of $\G(\gg_x)$ which consists of the equivalence
classes $[a]\in\G(\gg_x)$ of $\gg_x$-paths with the property
that, as an $\A$-path, $a$ is equivalent to the trivial $\A$-path.
\enddemo

The precise relation is as follows:

\proclaim{Lemma}
\label{N-versus-tilde-N}
For any Lie algebroid $A${\rm ,} and any $x\in M${\rm ,} $\tilde{N}_{x}(\A)$
is a subgroup of $\G(\gg_x)$ contained in the center
$Z(\G(\gg_x))${\rm ,} and its intersection with the connected component
$Z(\G(\gg_x))^{0}$ of the center is isomorphic to $N_{x}(\A)$.
\endproclaim 

\demo{Proof} 
Given $g\in \tilde{N}_{x}(\A)\subset \G(\gg_x)$ represented by a
$\gg_{x}$-path $a$, Proposition \ref{hol-hom} implies that parallel
transport $T_{a}:\gg_x\to \gg_x$ along $a$ is the identity. On the
other hand, since $a$ sits inside $\gg_x$, it is easy to see that
$T_{a}= {\rm ad}_{g}$, the adjoint action by the element $g\in \G(\gg_x)$
represented by $a$. This shows that $g\in Z(\G(\gg_x))$. The last part
follows from the fact that the exponential map induces an isomorphism
$\exp: Z(\gg_{x})\to Z(\G(\gg_x))^{0}$ (cf., e.g., $1.14.3$ in
\cite{DuKo}), and $N_{x}(\A)= \exp^{-1}(\tilde{N}_{x}(\A))$.
\enddemo

Since the group $\tilde{N}_{x}(A)$ is always countable (see next
section), we obtain:

\proclaim{{C}orollary}  For any Lie algebroid $A${\rm ,} and any $x\in M${\rm ,} the following are
equivalent\/{\rm :}
\begin{itemize}
\item[{\rm (i)}] $\tilde{N}_{x}(A)$ is closed\/{\rm ;}\/
\item[{\rm (ii)}] $\tilde{N}_{x}(A)$ is discrete\/{\rm ;}\/
\item[{\rm (iii)}] $N_{x}(A)$ is closed\/{\rm ;}\/
\item[{\rm (iv)}] $N_{x}(A)$ is discrete.
\end{itemize}
\endproclaim 

We remark that a special case of our main theorem shows that the
previous assertions are in fact equivalent to the integrability of
$A|_{L_x}$, the restriction of $\A$ to the leaf through $x$.

\demo{{\rm 3.2.} A second\/{\rm -}\/order monodromy map}
Let $L\subset M$ denote the leaf\break through~$x$. We define a homomorphism
$\partial:\pi_{2}(L,x)\to \G(\gg_x)$ with image precisely the
group $\tilde{N}_x(\A)$. This second-order monodromy map relates the
topology of the leaf through $x$ with the Lie algebra $\gg_x$.

Let $[\gamma]\in \pi_{2}(L,x)$ be represented by a smooth path
$\gamma: I\times I\to L$ which maps the boundary into $x$. We choose a
morphism of algebroids
$$  adt+ bd\epsilon : TI\times TI\to A_{L} $$ 
(i.e.\ $(a, b)$ satisfies equation (\ref{diffeq})) which lifts
$d\gamma: TI\times TI\to TL$ via the anchor, and such that $a(0,t)$,
$b(\epsilon,0)$, and $b(\epsilon,1)$ vanish. This is always possible:
for example, we can put $b(\epsilon, t)= \sigma
(\frac{d}{d\epsilon}\gamma(\epsilon, t))$ where $\sigma: TL\to A_{L}$
is a splitting of the anchor map, and take $a$ to be the unique
solution of the differential equation (\ref{diffeq}) with the initial
conditions $a(0, t)= 0$. Since $\gamma$ is constant on the boundary,
$a_{1}= a(1, -)$ stays inside the Lie algebra $\gg_{x}$, i.e.\ defines
a $\gg_{x}$-path
$$   a_{1}: I\to \gg_{x} \ .$$ 
Its integration (cf.\ \cite{DuKo}, or our Proposition \ref{A-G-paths}
applied to the Lie algebra $\gg_{x}$) defines a path in $\G(\gg_x)$,
and its end-point is denoted by $\partial(\gamma)$.
\enddemo

\proclaim{Proposition}\label{monodromy-seq}
The element $\partial(\gamma)\in \G(\gg_x)$ does not depend on the
auxiliary choices we made{\rm ,} and only depends on the homotopy class of
$\gamma$. Moreover{\rm ,} the resulting map \advance\eqcount by 3
\begin{equation}
\label{monodromy-map}
\partial: \pi_{2}(L,x)\to \G(\gg_x)
\end{equation}
is a morphism of groups and its image is precisely $\tilde{N}_{x}(\A)$.
\endproclaim 

Notice the similarity between the construction of $\partial$ and the
construction of the boundary map of the homotopy long exact sequence
of a fibration: if we view $0\to \gg_{L}\to A_{L}\to TL\to 0$ as
analogous to a fibration, the first few terms of the associated long
exact sequence will be
$$  \dots \to \pi_{2}(L,x)\stackrel{\partial}{\to} \G(\gg_x)
\to \G(A)_x\to \pi_1(L,x) .$$ 

The exactness at $\G(\gg_x)$ is precisely the last statement of the
proposition. We leave to the reader the (easy) check of exactness at
$\G(A)_x$. 

\demo{Proof of Proposition {\rm \ref{monodromy-seq}}}
{From} the definitions it is clear that
$\hbox{Im}~\partial=\tilde{N}_{x}(\A)$ so all we have to check is that
$\partial$ is well defined. For that we assume that
$$  \gamma^{i}= \gamma^{i}(\epsilon, t): I\times I\to L, \ \ i\in \{0, 1\} $$ 
are homotopic relative to the boundary, and that
$$  a^{i}dt+ b^id\epsilon: TI\times TI\to A_{L}.\ \ i\in \{0, 1\} $$ 
are lifts of $d\gamma^{i}$ as above. We prove that the paths
$a^{i}(1, t)$ ($i\in \{0, 1\}$) are homotopic as $\gg_x$-paths.

By hypothesis, there is a homotopy $\gamma^{u}=
\gamma^{u}(\epsilon, t)$ ($u\in I$) between $\gamma^{0}$ and
$\gamma^{1}$. We choose a family $b^{u}(\epsilon, t)$ joining
$b^{0}$ and $b^{1}$, such that $\#(b^{u}(\epsilon, t))=
\frac{d\gamma^u}{d\epsilon}$ and $b^{u}(\epsilon, 0)= b^{u}(\epsilon,
1)= 0$. We also choose a family of sections $\eta$ depending on $u,
\epsilon, t$ such that
$$ 
\eta^{u}(\epsilon, t, \gamma^{u}(\epsilon, t))= b^{u}(\epsilon, t),\hbox{
with }
\eta= 0 \hbox{ when } t= 0, 1.
$$ 
As in the proof of Proposition \ref{equivalence}, let $\xi$ and
$\theta$ be the solutions of
$$ \left\{
\begin{array}{l}
\frac{d\xi}{d\epsilon}-\frac{d\eta}{dt}= [\xi, \eta],\hbox{ with }
\xi=0 \hbox{ when } \epsilon = 0, 1,\\  \\
\frac{d\theta}{d\epsilon}-\frac{d\eta}{du}=[\theta, \eta],\hbox{ with }
\theta=0 \hbox{ when } \epsilon = 0, 1.
\end{array}
\right.$$ 
Setting $u= 0,1$ we get
$$  a^{i}(\epsilon, t)=
\xi^{i}(\epsilon, t, \gamma^{i}(\epsilon, t)), \quad i=0, 1.$$ 
On the other hand, setting $t= 0, 1$ we get
$$ 
\theta= 0\hbox{ when } t= 0, 1.
$$ 
A brief computation shows that
$\phi\equiv\frac{d\xi}{du}-\frac{d\theta}{dt}-[\xi,\theta]$
satisfies
$$ 
\frac{d\phi}{d\epsilon}=[\phi, \eta],
$$ 
and since $\phi =0$  when $\epsilon= 0$, it follows that
$$  \frac{d\xi}{du}- \frac{d\theta}{dt}= [\xi, \theta] .$$ 
If in this relation we choose $\epsilon=1$, and use
$\theta^u(1,t)=0$ when $t=0,1$, we conclude that $a^{i}(1, t)=
\xi^i(1,t,\gamma^{i}(1, t))$, $i= 0,1$, are equivalent as $\gg_x$-paths. \phantom{hunger}
\enddemo

3.3. {\it Computing the monodromy}. Let us indicate briefly how the monodromy groups (Definition
\ref{N-groups} or, alternatively, Definition \ref{tilde-N-groups}),
can be explicitly computed in many examples. We consider the short
exact sequence
$$  0\to \gg_{L}\to \A_{L}\stackrel{\#}{\to} TL\to 0 $$ 
and a linear splitting $\sigma: TL\to \A_{L}$ of $\#$. The
{\it curvature} of $\sigma$ is the element $\Omega_{\sigma}\in
\Omega^{2}(L;\gg_{L})$ defined by:
$$  \Omega_{\sigma}(X, Y)\equiv \sigma ([X, Y])- [\sigma(X), \sigma(Y)]\ .$$ 
In favorable cases, the computation of monodromy can be reduced to the
following:

\proclaim{Lemma}
\label{compute}
If there is a splitting $\sigma$ with the property that its
curvature $\Omega_{\sigma}$ is $Z(\gg_L)$\/{\rm -}\/valued{\rm ,} then
$$  N_{x}(A)= \{ \int_{\gamma} \Omega_{\sigma}:
[\gamma]\in \pi_{2}(L, x)\}\subset Z(\gg_x)$$ 
for all $x\in L$.
\endproclaim 

Before we give a proof some explanations are in order.

First of all, $Z(\gg_{L})$ is canonically a flat vector bundle over
$L$. The corresponding flat connection can be expressed with the
help of the splitting $\sigma$ as
$$  \nabla_{X}\alpha= [\sigma(X), \alpha ],$$ 
and it is easy to see that the definition does not depend on
$\sigma$. In this way $\Omega_{\sigma}$ appears as a $2$-cohomology
class with coefficients in the local system defined by $Z(\gg_{L})$
over $L$, and then the integration is just the usual pairing
between cohomology and homotopy. In practice one can always avoid
working with local coefficients: if $Z(\gg_{L})$ is not already
trivial as a vector bundle, one can achieve this by pulling back to
the universal cover of $L$ (where parallel transport with respect
to the flat connection gives the desired trivialization).

Second, we should specify what we mean by integrating forms with coefficients in a local system.  Assume $\omega\in
\Omega^2(M;E)$ is a 2-form with coefficients in some flat vector bundle $E$.  Integrating $\omega$ over a 2-cycle
$\gamma:{\Bbb S}^2\to M$ means (i) taking the pull-back $\gamma^\ast\omega\in \Omega^2({\Bbb S}^2; \gamma^\ast
E)$, and (ii) integrate $\gamma^\ast \omega$ over ${\Bbb S}^2$.  Here $\gamma^\ast E$ should be viewed as a flat vector
bundle of $\Bbb S^2$ for the pull-back connection.  Notice that the connection enters the integration part, and this
matters for the integration to be invariant under homotopy.


\demo{Proof  of Lemma {\rm \ref{compute}}}
We may assume that $L= M$, i.e.\ $\A$ is transitive. In agreement with the
comments above, we also assume for simplicity that $Z({\frak g})$ is
trivial as a vector bundle $({\frak g} \equiv {\frak g}_L)$.
The formula above defines a connection $\nabla^{\sigma}$ on the entire
${\frak g}$. We use $\sigma$ to identify $\A$ with $TM\oplus
{\frak g}$
so the bracket becomes
$$  [(X, v), (Y, w)]= ([X, Y], [v, w]+ \nabla_{X}^{\sigma}(w)-
\nabla_{Y}^{\sigma}(v)- \Omega_{\sigma}(X, Y)) .$$ 
We choose a connection $\nabla^{M}$ on $M$, and we consider the
connection $\nabla= (\nabla^{M}, \nabla^{\sigma})$ on $\A$. Note that
$$  T_{\nabla}((X, v), (Y, w))= (T_{\nabla^{M}}(X, Y), \Omega_{\sigma}(X, Y)-
[v, w]) $$ 
for all $X, Y\in TM$, $v, w\in {\frak g}$. This shows that two
$\A$-paths $a$ and $b$ as in Proposition \ref{equivalence} will be of
the form $a=(\frac{d\gamma}{dt}, \phi)$, $b= (\frac{d\gamma}{d\epsilon},
\psi)$
where $\phi,\psi$ are paths in ${\frak g}$ satisfying
$$  \partial_{t}\psi- \partial_{\epsilon}\phi=
\Omega_{\sigma}(\frac{d\gamma}{dt}, \frac{d\gamma}{d\epsilon})- [\phi,
\psi].$$ 
Now we only have to apply the definition of $\partial$: Given
$[\gamma]\in \pi_{2}(M, x)$, we choose the lift $adt+ bd\epsilon$ of
$d\gamma$ 
with $\psi= 0$ and 
$$  \phi = -\int_{0}^{\eps}
\Omega_{\sigma}(\frac{d\gamma}{dt}, \frac{d\gamma}{d\epsilon}) .$$ 
Then $\phi$ takes values in $Z({\frak g}_{x})$, and we obtain
$\partial [\gamma]=[\int_{\gamma}\Omega_{\sigma}]$.
\enddemo

\numbereddemo{{E}xample}
\label{ex:central:extension}
Recall (e.g.\ \cite{Mack1}) that any closed two-form $\omega\in
\Omega^2(M)$ induces an algebroid $\A_{\omega}=TM\oplus{\Bbb L}$, where ${\Bbb L}$ is the trivial line bundle, 
with anchor $(X,\lambda)\mapsto X$ and Lie bracket
$$  [(X,f),(Y,g)]=([X, Y],X(g)-Y(f)+\omega(X,Y)).$$ 
Using the obvious splitting of $\A$, Lemma 3.3 tells
us that
$$  N_{x}(\A_{\omega})= \set{ \int_{\gamma} \omega:
     [\gamma]\in \pi_{2}(M, x)}\subset {\Bbb R} $$ 
is the group of periods of $\omega$. Other examples will be discussed
in the next sections.
\enddemo

3.4. {\it Measuring the monodromy}.
 In order to measure the size of the monodromy groups $N_x(A)$,
we fix some norm on the Lie algebroid $\A$ and for $x\in M$ we set
$$  r(x)\equiv d(0, N_{x}(\A)-\{0\}), $$ 
where we adopt the convention that $d(0,\emptyset)=+\infty$.

When $x$ varies on a leaf $L$ this function varies continuously,
since the norm on $A$ is assumed to vary continuously and the
groups $N_{x}(A)$ are all isomorphic for $x\in L$. On the other
hand, when $x$ varies in a transverse direction the behavior of
$r(x)$ is far from being continuous as illustrated by the following
examples:

\numbereddemo{{E}xample}
\label{example:basic}
We take for $A$ the trivial $3$-dimensional vector bundle over $M=\Rr^3$,
with basis
$\set{e_1,e_2,e_3}$. The Lie bracket on $A$ is defined by
\begin{eqnarray*}
\left[e_2,e_3\right]&=&a e_1+ b x^1 \bar{n},\\
\left[e_3,e_1\right]&=&a e_2+ b x^2 \bar{n},\\
\left[e_1,e_2\right]&=&a e_3+ b x^3 \bar{n}
\end{eqnarray*}
where $\bar{n}= \sum_i x^ie_i$ is a central element, and
depends on two (arbitrary) smooth functions $a$ and $b$ of the
radius $R$, with $a(R)>0$ whenever $R>0$. The
anchor is given by
$$  \#(e_i)=a v_{i}, \qquad i=1,2,3 $$ 
where $v^i$ is the infinitesimal generator of a rotation about the
$i$-axis:
$$ 
v_1=x^3\frac{\partial}{\partial x^2}-x^2\frac{\partial}{\partial x^3},\quad
v_2=x^1\frac{\partial}{\partial x^3}-x^3\frac{\partial}{\partial x^1},\quad
v_3=x^2\frac{\partial}{\partial x^1}-x^1\frac{\partial}{\partial x^2}.$$ 
The leaves of the foliation induced on ${\Bbb R}^3$ are the
spheres $S_{R}^{2}$ centered at the origin, and the origin is the
only singular point.

We now compute the function $r$ using the obvious metric on $\A$.
We restrict to a leaf $S_{R}^{2}$ with $R>0$, and as splitting of $\#$ we
choose the map defined by
$$  \sigma(v_i)= \frac{1}{a}\left(e_i-\frac{x^i}{R^2}\bar{n}\right).$$ 
Then we obtain the center-valued 2-form (cf.~\S 3.3)
$$  \Omega_{\sigma}= \frac{bR^2- a}{a^2R^4}\omega \bar{n}$$ 
where $\omega= x^1 dx^2\wedge dx^3+ x^2dx^3\wedge dx^1+ x^3dx^1\wedge dx^2$.
Since $\int_{S^2_R}\omega=4\pi R^3$ it follows that
$$   N(\A)\simeq 4\pi \frac{bR^2- a}{a^2 R} \Zz\bar{n}\subset \Rr\bar{n}.$$ 
This shows that
$$ 
r(x,y,z)= \left\{
\begin{array}{ll}
+\infty \qquad&\  \hbox{if}\ R=0\ \hbox{ or }  a=b R^2,\\
\\
4\pi \frac{bR^2-a}{a^2} \qquad&\  \hbox{otherwise.}
\end{array}
\right.
$$ 
So the monodromy might vary in a nontrivial fashion, even nearby
regular leaves.
\enddemo

In the previous example the function $r$ is {\it not} upper
semi-continuous. In the next example we show that $r$, in general, is
{\it not} lower semi-continuous. This example also shows that, even if the
anchor is injective in a ``large set'', one has no control on the way
the monodromy groups vary.

\numbereddemo{{E}xample}
We consider a variation of Example \ref{example:basic}, so we use the
same notation.  We let $M= {\Bbb S}^{2}\times {\Bbb H}$, where
${\Bbb H}$ denotes the quaternions. The Lie algebroid $\pi:A\to M$
is trivial as a vector bundle, has rank $3$, and relative to a basis
of sections $\set{e_1, e_2, e_3}$ the Lie bracket is defined by $[e_1,
e_2]= e_3$ and cyclic permutations. To define the anchor, we let $v_1,
v_2, v_3$ be the vector fields on ${\Bbb S}^2$ obtain by restricting
the infinitesimal generators of rotations, and we let $w_1, w_2, w_3$
be the vector fields on ${\Bbb H}$ corresponding to multiplication
by $\vec{i},\vec{j}, \vec{k}$. The anchor of the algebroid is then
defined by setting $\# e_i\equiv(v_i, w_i)$, $i=1,2,3$. For this Lie
algebroid one has:
\begin{itemize}
\item the anchor is injective on a dense open set;
\item there is exactly one singular leaf, namely the sphere
${\Bbb S}^2\times \{ 0\}$.
\end{itemize}
Now observe that the monodromy above the singular leaf is nontrivial,
since the restriction of $A$ to this singular leaf is the central extension
algebroid $T{\Bbb S}^2\oplus {\Bbb L}$ defined by the area form on
${\Bbb S}^2$. For the function $r$ we have again:
$$ 
r(x)= \left\{
\begin{array}{ll}
r_0 \qquad&\  \hbox{if}\ x\in {\Bbb S}^2\times \{ 0\}\\
\\
+\infty \qquad&\  \hbox{otherwise.}
\end{array}
\right.
$$ 
Note that in this case $r_0>0$. This is no accident, as we will
see later in Section~5.2.5.
\enddemo

Finally, we give an example where a splitting as in Lemma
\ref{compute} does not exist, and which illustrates how the
groups $N_{x}$ and $\tilde{N}_{x}$ can differ.

\numbereddemo{{E}xample}
Let $\gg$ be any Lie algebra and $G$ a simply connected Lie group with
Lie algebra $\gg$. If we consider the Poisson manifold $M=\gg^*$ with the
Kirillov-Kostant Poisson bracket, and let $A=T^*\gg^*$ (for the cotangent Lie algebroid of a Poisson manifold; see [4]).
We always obtain an integrable Lie algebroid: a source-simply connected Lie groupoid
integrating $A$ is $\G=T^*G\simeq G\times \gg^*$ with source and target maps
$$  \s(g,\xi)=\xi,\qquad \t(g,\xi)=\Ad^*g\cdot\xi,$$ 
and with multiplication
$(g_1,\xi_1)\cdot(g_2\cdot\xi_2)=(g_1g_2,\xi_2)$, wherever defined.

For a specific example take $M={\frak su}^*(3)$. The symplectic
leaves (i.e., the co-adjoint orbits) are isospectral sets, and so we can
understand them by looking at their point of intersection with the
diagonal matrices with imaginary eigenvalues. There are orbits of
dimension 6 (distinct eigenvalues), dimension~4 (two equal
eigenvalues) and the origin (all eigenvalues equal). Let us take for example
the (singular) orbit $L$ through
$$  x=\left(\begin{array}{ccc}
i\lambda &0 &0\\
0&i\lambda&0\\
0&0&-2i\lambda
\end{array}\right).$$ 
Then we find its isotropy subalgebra to be
$$  
\gg_x=\set{
\left(
\begin{tabular}{c|c}
$\begin{array}{ccc}
& &\\
&X&\\
& &
\end{array}$
& 0\\ \cline{1-2}
0& $-\tr X$
\end{tabular}
\right):\ X\in {\frak u}(2)}
$$ 
and so we see that the simply connected Lie group integrating the Lie
algebra $\gg_x$ is $\G(\gg_x)=\Rr\times {\rm SU}(2)$. On the other hand, the
isotropy group $\G_x$ is given by
\begin{eqnarray*} 
\G_x&=&\set{g\in {\rm SU}(3): \Ad^*g\cdot x=x}\\
&=&\set{
\left(
\begin{tabular}{c|c}
$\begin{array}{ccc}
& &\\
&g&\\
& &
\end{array}$
& 0\\ \cline{1-2}
0& $\det g^{-1}$
\end{tabular}
\right):\ g\in {\rm U}(2)}.
\end{eqnarray*}
We conclude that the orbit $L$ is diffeomorphic to
${\rm SU}(3)/{\rm U}(2)={\rm CP}(2)$. In fact, one can show that it is symplectomorphic
to ${\rm CP}(2)$ with its standard symplectic structure (see \cite[Example 3.4.5]{Fer2}). Also, we see that the long exact
sequence
$$  \dots \to \pi_{2}({\rm CP}(2),x)\stackrel{\partial}{\to} \G(\gg_x)
\to \G(T^*M)_x\to \pi_1({\rm CP}(2),x),$$ 
reduces to:
$$  \dots \to \Zz\stackrel{\partial}{\to} \Rr\times {\rm SU}(2)
\stackrel{\rho}{\to} {\rm U}(2)\to \set{1},$$ 
where $\rho(\theta,A)=e^{i\theta}A$. We conclude that $\partial n=(\pi
n,(-1)^nI)$, so that $\partial$ takes values in the center $Z(\Rr\times
{\rm SU}(2))=\Rr\times \set{\pm I}$, and
$$  \tilde{N}_x=\hbox{Im~}\partial=2\Zz\times\set{\pm I},\qquad N_x=2\Zz.$$ 
Since these two groups are distinct, there can be no splitting as in Lemma
\ref{compute}. Another argument is that such a splitting would define
a flat connection on the conormal bundle $\nu^*(L)=\Ker\#|_L$, and since
$L={\rm CP}(2)$ is 1-connected, it would follow that the conormal bundle would be
a trivial bundle. This is not possible. In fact, the total
Stiefel-Whitney class of ${\rm CP}(2)$ is nontrivial, and hence so is the
Stiefel-Whitney class of the normal bundle, when we embed ${\rm CP}(2)$ in
any euclidean space. So the conormal bundle cannot be trivial.
\enddemo

\vglue-8pt
%
\section{Obstructions to integrability} %
\label{section:integrability}           %
\vglue-4pt

In this section we first state our main result which gives the
obstructions to integrability, and give a few examples. We then
give another description of the Weinstein groupoid which is more
suitable for proving the theorem.

\demo{{\rm 4.1.} The main theorem}            
Let $\A$ be a Lie algebroid over $M$. Using the notations introduced
above, our main result is the following:
\enddemo

\proclaim{Theorem}
\label{main thm}
A Lie algebroid $\A$ over $M$ is integrable if and only if\/{\rm :}
\begin{itemize}
\item[{\rm (i)}] Longitudinal obstruction: $N_{x}(\A) \subset \A_{x}$ is discrete
{\rm (}\/i.e.{\rm ,} $r(x)\neq 0${\rm ),}
\item[{\rm (ii)}] Transverse obstruction\/{\rm :} $\liminf_{y\to x}r(y)>0${\rm ,}
\end{itemize}
for all $x\in M$.
\endproclaim 

The next examples illustrate this result and show that these two
obstructions are independent.

\numbereddemo{{E}xample}
\label{ex:first:obstruction}
In this example, nonintegrability is forced by the first obstruction.
We simply take the central extension Lie algebroid
$\A_{\omega}=TM\oplus{\Bbb L}$ associated with a closed 2-form on
$M$ with a noncyclic group of periods (cf.~Example
\ref{ex:central:extension}). Then $r(x)=0$ so the first obstruction
ensures us that $\A_{\omega}$ is nonintegrable. We point out that
this is a well-known counter-example to integrability (cf.~e.g.~
\cite[p.~118]{CaWe}) which is usually approached through the theory of
transversally parallelizable foliations (see also \S
5.3
below).
\enddemo

\numbereddemo{{E}xample}
\label{ex:second:obstruction}
Let us give an example of a regular Lie algebroid with trivial first
obstruction, while the second one is not. Take $\F$ to be the trivial
foliation of $M=N\times T$ with leaves $N\times\{t\}$, $t\in T$. Also
we choose a closed 2-form $\omega$ on $N$ with cyclic group of periods
and we set $\omega_t=\phi(t)\omega$, where $\phi$ is some smooth
function on $T$. Since the pull-back of $\omega_t$ to any leave is
closed, we obtain the central extension Lie algebroid
$A_{\omega_t}=\F\oplus {\Bbb L}$, as in Example
\ref{ex:central:extension}, the leaves of which are the leaves of
$\F$.

The ``first obstruction'' is satisfied for all leaves, but clearly the
``second obstruction'' is not satisfied at the points $t_0\in T$ with
the property that\break $\phi(t_0)= 0$ and $\phi$ is not locally constant at
$t_0$.
\enddemo

\numbereddemo{{E}xample}
\label{example:second:obstruction:2}
Consider the Lie algebroid $A$ over $\Rr^3$ discussed in
Example~\ref{example:basic}. Then $\A$ satisfies the first obstruction, but
it does not satisfy the second obstruction at points where $aR^2-b$
vanishes (without vanishing identically in some neighborhood of the
point) and also at the origin if $\liminf_{R\to 0}\frac{b R^2-
a}{a^2}=0$.

For example, if we choose $a= R^2$, $b= R^3+1$, the resulting Lie
algebroid $A$ over $\Rr^3$ has the following two properties:
\begin{itemize}
\item[(a)] Its restriction to ${\Bbb R}^{3}-{0}$ is integrable;
\item[(b)] Its restriction to any disk around the origin is not integrable
(because of the second obstruction at $x=0$).
\end{itemize}

\enddemo

\numbereddemo{{E}xample}
\label{example:Weinstein}
Let us explain Weinstein's example of a nonintegrable regular Poisson
manifold given in \cite{Wein} (see also \cite[\S 6]{CaFe}).  He
takes $M={\Bbb R}^3-\{0\}\cong {\frak su}(2)^*-\{0\}$ with the
Kirillov-Poisson structure scaled by a function $f(R)$ depending
on the radius.  The associated cotangent Lie algebroid is in fact
$T^*M=A|_{{\Bbb R}^3-\{0\}}$, where $A$ is the Lie algebroid of
Examples \ref{example:basic} and \ref{example:second:obstruction:2}
with $a=f$, $b= \frac{1}{R}f'$. Its integrability is then controlled by
$$ r(R)=4\pi \frac{Rf'-f}{f^2}=-A'(R),$$ 
where $A(R)=\frac{4\pi R}{f}$ is the symplectic area.
\enddemo

We refer to Section \ref{Examples/Applications} for various integrability
criteria that can be deduced from the theorem, including
all criteria that have appeared before in the literature.

\demo{{\rm 4.2.} The Weinstein groupoid as a leaf space} Before we can proceed with the proof of our main result, we
need a better control on the equivalence relation defining the Weinstein
groupoid $\G(\A)$. In this section we will show that $\G(\A)$ is
the leaf space of a foliation $\F(\A)$ on $P(\A)$, of
finite codimension, the leaves of which are precisely the equivalence
classes of the homotopy relation $\sim$ of $\A$-paths.

As before, $\A$ is a fixed Lie algebroid over $M$.  We will use the
following notations when working in local coordinates: we let
$x=(x^{1},\dots,x^{n})$ denote local coordinates on $M$, and we denote
by $\set{e_{1},\dots,e_ {k}}$ a (local) basis of $A$ over this
chart. The anchor and the bracket of $A$ decompose as
$$  \#e_{p}= \sum_i b^{i}_{p}\frac{\partial}{\partial x^{i}},\qquad
   [e_p, e_q]= \sum_{r} c_{pq}^{r} e_{r},$$ 
and an $\A$-path $a$ can be written as $a(t)= \sum_{p} a^{p}(t)e_p$.

Let us first describe the smooth structure on $P(\A)$. We consider the
larger space $\tilde{P}(\A)$ of all $C^1$-curves $a: I\to \A$ with
base path $\gamma= \pi\circ a$ of class $C^2$. It has an obvious
structure of Banach manifold, and its tangent space
$T_{a}(\tilde{P}(\A))$ consists of curves $U:I\to TA$ such that $U(t)\in
T_{a(t)}A$. Using a $TM$-connection $\nabla$ on $\A$, such curves can
be viewed as pairs $(u, \phi)$ formed by a curve $u:I\to A$ over
$\gamma$ and a curve $\phi:I\to TM$ over $\gamma$ (namely, the
vertical and horizontal component of $U$).

\proclaim{Lemma}\label{P=Banach} $P(\A)$ is a {\rm (}\/Banach\/{\rm )} submanifold of
$\tilde{P}(\A)$. Moreover{\rm ,} given a connection $\nabla$ on $\A${\rm ,} the
tangent space $T_{a}P(\A)$ consists of the paths $U= (u, \phi)$
such that
$$ 
\#u= \overline{\nabla}_{a}\phi.
$$ 
\endproclaim 
\vglue-9pt
{\it Proof}.
We consider the smooth map $F:\tilde{P}(\A)\to \tilde{P}(TM)$ given
by
$$  \ F(a)=\#a-\frac{d}{dt}\pi\circ a. $$ 
Clearly $P(\A)= F^{-1}(Q)$, where $Q$ is the subspace of
$\tilde{P}(TM)$ consisting of zero paths. Fix $a\in P(\A)$, with
base path $\gamma=\pi\circ a$, and let us compute the image of $U=
(u,\phi)\in T_{a}\tilde{P}(\A)$ under the differential
$$  (dF)_{a}: T_{a}\tilde{P}(\A)\to  T_{0_{\gamma}}\tilde{P}(TM) \ .$$ 
The result will be some path $t\mapsto (dF)_{a}\cdot U(t)\in
T_{0_{\gamma(t)}}TM$; hence, using the canonical splitting
$T_{0_x}TM\cong T_{x}M\oplus T_{x}M$, it will have a horizontal and
vertical component. We claim that for any connection $\nabla$, if
$(u, \phi)$ are the components of $U$, then
$$ 
((dF)_{a}\cdot U)^{\hbox{hor}}= \phi,\quad
((dF)_{a}\cdot U)^{\hbox{ver}}= \#u
-\overline{\nabla}_{a}\phi.
$$ 
Note that this immediately implies that $F$ is transverse to $Q$,
so the assertion of the proposition follows. Since this
decomposition is independent of the connection $\nabla$ and it is
local (we can look at restrictions of $a$ to smaller intervals), we
may assume that we are in local coordinates, and that $\nabla$ is
the standard flat connection. We now use the notation above, and
we denote by $\frac{\partial}{\partial x^{i}}$ the horizontal basis
of $T_{0_x}TM$, and by $\frac{\delta}{\delta x^{i}}$ the vertical
basis. A simple computation shows that the horizontal component of
$(dF)_{a}(u,\phi)$ is $\sum \phi^{i}\frac{\partial}{\partial
x^{i}}$, while its vertical component is
$$  \sum_{j} \left(
-\dot{\phi}^{j}(t)+ \sum_{p} u^{p}(t) b_{p}^{j}(\gamma(t))
+\sum_{p, i} a^{p}(t)\phi^{i}(t)\frac{\partial b_{p}^{j}}{\partial
x^i}(\gamma(t))
\right) \frac{\delta}{\delta x^{j}}.$$ 
That this is precisely $\#u - \overline{\nabla}_{a}\phi$ immediately
follows by computing
\begin{eqnarray*}
\overline{\nabla}_{e_{p}}\frac{\partial}{\partial x^{i}}
    & =&  \# \nabla_{\frac{\partial}{\partial x^{i}}}e_{p}-
         \left[\frac{\partial}{\partial x^{i}}, \#e_p\right]\\
    & = &-\sum_{j} \frac{\partial b_{p}^{j}}{\partial x^i}
         \frac{\partial}{\partial x^{j}}.\\
\noalign{\vskip-36pt}
\end{eqnarray*}
\phantom{ha}\hfill\qed
\vglue12pt

We now construct an involutive sub-bundle $\F(\A)$ of
$TP(\A)$, i.e.\ a foliation on $P(\A)$. Let us fix a connection
$\nabla$ on $\A$, and let $a$ be an $\A$-path with base path
$\gamma$. We denote by $\Pa_{0, \gamma}(\A)$ the space
of all $C^2$-paths $b: I\to A$ such that $b(t)\in A_{\gamma(t)}$ and
$b(0)= b(1)= 0$. For any such $b$, we have a tangent vector
$X_{b, a}\in T_{a}P(\A)$ with components $(u, \phi)$ relative
to the connection $\nabla$ given by
$$ 
u= \overline{\nabla}_{a}b,\quad \phi= \#b.
$$ 
Lemma \ref{P=Banach} shows that these are indeed tangent to $P(\A)$, and we
set:
$$  \F_{a}(\A)\equiv \set{ X_{b, a}\in T_{a}P(\A): b\in  \Pa_{0,
\gamma}(\A)}.$$  

Some geometric insight to this sub-bundle can be obtain by
considering the Lie algebra of time-depending sections of $A$ vanishing at
the end-points:
$$  P_{0}\Gamma(\A)= \set{ I\ni t\mapsto \eta_{t}\in \Gamma(\A): \eta_{0}
                  = \eta_{1}=0, \eta\hbox{ is of class } C^2\hbox{ in } t}.
$$ 
For any such section $\eta$ we consider the induced path $b(t)=
\eta(t, \gamma(t))$ and put $X_{\eta, a}\equiv X_{b, a}$. The
resulting map
$$   P_{0}\Gamma(\A) \to \X(P(\A)),\quad \eta\mapsto X_{\eta} $$ 
is an action of the Lie algebra $P_{0}\Gamma(\A)$ on $P(\A)$.

 
We now show that the foliation $\F(A)$ is in fact the same as the
partition of $P(A)$ into equivalent classes \pagebreak of $A$-paths:

\proclaim{Proposition}
\label{G-as-fol}
For a Lie algebroid $A${\rm :}
\begin{itemize}
\item[{\rm (i)}] The spaces $\F_{a}(\A)$ do not depend on the choice of the
connection $\nabla$. More precisely{\rm ,} for any $\eta\in
P_{0}\Gamma(\A)${\rm ,}
$$ 
X_{\eta, a}(t)= {\left. \frac{d}{d\epsilon}\right|}_{\epsilon= 0}
\phi_{\eta_{t}}^{\epsilon, 0} a(t)+ \frac{d\eta_{t}}{dt} (\gamma(t)).
$$ 
\item[{\rm (ii)}] $\F(\A)$ is a foliation on $P(\A)$ of finite codimension
equal to $n+k$ where $n=\dim M$ and $k=\rank A$.
\item[{\rm (iii)}] Two $\A$\/{\rm -}\/paths are equivalent {\rm (}\/homotopic\/{\rm )} if and only if
they are in the same leaf of $\F(\A)$.
\item[{\rm (iv)}] For any {\rm (}\/local\/{\rm )} connection $\nabla$ on $\A${\rm ,} the
exponential map $\Exp_{\nabla}: \A\to P(\A)$ is transverse to
$\F(\A)$.
\end{itemize}

\endproclaim 

\demo{Proof}  
We first assume that $\eta_t$ is an element of  
$P_{0}\Gamma(\A)$  and we will see
that it induces a vector field $X_{\eta}$ on $P(\A)$ tangent to
$\F(\A)$ the flow of which preserves the equivalence of paths. This is
only a reformulation of Proposition \ref{equivalence}. Hence, let
$a_{0}\in P(\A)$ with base path $\gamma_{0}$, and let $\xi_{0}$ be a
time-dependent section of $\A$ such that $\xi_{0}(t, \gamma_{0}(t))=
a_{0}(t)$.  We denote by $\xi$ the solution of (\ref{eq.equiv}) with
the initial condition $\xi(0, t)= \xi_{0}(t)$. Then, as in the proof
of Proposition~\ref{equivalence},
$$  \xi(\epsilon, t)= \int_{0}^{\epsilon} (\phi_{\eta_{t}}
                     ^{\epsilon', \epsilon})^{*}
                     \frac{d\eta}{dt} d\epsilon'
                   + (\phi_{\eta_{t}}^{0, \epsilon})^{*}\xi_{0} .$$ 
Now consider the base path
$$  \gamma_{\epsilon}(t)= \Phi_{\#\eta_{t}}^{\epsilon, 0}\gamma_{0}(t) $$ 
and the paths above it
$$  a_{\epsilon}(t)= \xi_{\epsilon}(t, \gamma_{\epsilon}(t)),\
b_{\epsilon}(t)= \eta (t, \gamma_{\epsilon}(t))\ .$$ 
We can view $\epsilon\mapsto a_{\epsilon}$ as a curve in $P(\A)$
starting at $a_{0}$, and defining a tangent vector
$$ 
{\left. \frac{d}{d\epsilon}\right|}_{\epsilon= 0} a_{\epsilon}
\in T_{a_{0}}P(\A) .
$$ 
Given some connection $\nabla$, Proposition \ref{equivalence} shows
that this tangent vector has vertical component
$$  \partial_{\epsilon}a= \partial_{t}b- T(a, b)= \overline{\nabla}_{a}(b) $$ 
at $\epsilon= 0$, while the horizontal component is
$$   {\left. \frac{d}{d\epsilon}\right|}_{\epsilon= 0}
\gamma_{\epsilon}(t)= \#b_{0}(t) .$$ 
In other words,
$$ {\left. \frac{d}{d\epsilon}\right|}_{\epsilon= 0}
a_{\epsilon}=X_{b_{0}, a_{0}}\in \F_{a_{0}}(\A).$$ 
On the other hand, the formula above describing $\xi$ shows that
$$  a_{\epsilon}(t)= \int_{0}^{\epsilon}
                    \phi_{\eta_{t}}^{\epsilon, \epsilon'}
                    \frac{d\eta_{t}}{dt}
                    (\phi_{\#\eta_{t}}^{\epsilon',
0}(\gamma_{0}(t))
                    d\epsilon'
                  + \phi_{\eta_{t}}^{\epsilon, 0}(\gamma_{0}(t)) .
$$ 
The derivative at $\epsilon= 0$ is precisely the expression given in (i) and
this also shows that (iii) holds.

To determine the codimension of $\F_{a}(\A)$ we note that given $(u,\phi)$
satisfying $\#u=\overline{\nabla}_{a}\phi$ (i.e., a vector tangent to
$P(\A)$)
and lying in $\F_{a}(\A)$, we have
\begin{itemize}
\item[(a)] $\phi(0)=0$;
\item[(b)] If we consider the solution $b$ of the equation
      $\overline{\nabla}_{a}(b)= u$ with initial condition
      $b(0)= 0$ (which can be expressed in terms of the parallel
      transport along $a$ with respect to $\overline{\nabla}$),
      we must have $b(1)= 0$.
\end{itemize}
Conversely, if (a) and (b) hold, we have that
$\overline{\nabla}_{a}(\#b-\phi)=0$ and $\#b-\phi$ vanishes at $t=
0$. It follows that $\phi=\#b$ and $u=\overline{\nabla}_{a}b$, so
$(u,\phi)$ is a tangent vector in $\F_{a}(\A)$. This shows that
$\codim \F=\dim M+\rank A$.

Finally, to prove (iv), we assume for simplicity that we are in
local coordinates and that $\nabla$ is the trivial flat connection
(this is actually all we will use for the proof of the main
theorem, and this in turn will imply the full statement of (iv)).
Also, all we need to show is that $\Exp_{\nabla}(A)$ is transverse
to $\F(\A)$ at any trivial $A$-path $a=O_{x}$ over $x\in M$. Now,
the equations for the geodesics show that if $(u,\phi)$ is a
tangent vector to $\Exp_{\nabla}(A)$ at $a$ then we must have:
$$  \dot{\phi}^i=b^i_p(x)u^p,\qquad \dot{u}^p=0.$$ 
Therefore, we see that:
$$  T_{a}\Exp_{\nabla}(\A)=\set{(u,\phi)\in T_{a}P(\A): u(t)=u_0,
\phi(t)=\phi_0+t\#u_0}.$$ 
Suppose that a tangent vector $(u,\phi)$ belongs to this $n+k$-dimensional space and is also tangent to $\F(\A)$. Then
(a) above implies that $\phi_0=0$, while (b) says that the solution $b$ of
$\frac{d b}{dt}=u_0$ with $b(0)=0$ satisfies $b(1)=0$.  Therefore,
we must have $\phi_0=0$ and $u_0=0$, so $(u,\phi)$ is the null
tangent vector. This shows that $\Exp_{\nabla}(A)$ is transverse to
$\F(\A)$ at $0_x$, for any $x$.
\enddemo

{\it Remark} 4.8.
In \cite{CaFe}, Cattaneo and Felder obtain the Weinstein groupoid
for the special case of Poisson manifolds by a Hamiltonian
reduction procedure. The Lie algebraic interpretation given above
for the foliation $\F(\A)$ shows that our construction of $\G(\A)$
for general $\A$ is also obtained by a kind of reduction procedure
for Lie algebroid actions. \pagebreak
  
\advance\theoremcount by 1

4.3. {\it Proof of the main theorem}.
In this section we prove Theorem \ref{main thm} (for notations, see
\S \ref{section:monodromy}).

To prove that both conditions are necessary, choose some connection
$\nabla$ on $\A$, and let $\overline{\Exp}_{\nabla}: A\to \G(\A)$
be the associated exponential map. Clearly the restriction of
$\overline{\Exp}^{\phantom{|}}_{\nabla}$ to $\gg_x$ is the composition of the
exponential map of $\gg_x$ with the obvious map $i: \G(\gg_x)\to
\G(A)_x$, which shows that $\overline{\Exp}_{\nabla}(v_x)= 1_x$
for all $v_x\in N_{x}(\A)$ in the domain of the exponential map. On
the other hand, if $\A$ is integrable, we know that
$\overline{\Exp}_{\nabla}$ will be a diffeomorphism on a small
neighborhood of $M$ on $\A$. Hence there must exist an open
$\U\subset \A$ such that $\U\cap N(\A)= M$, where $N(\A)=
\mathbold{\cup}_{x}N_{x}(\A)$. But this is obviously equivalent to the
conditions in the statement.

We now show that these conditions also guarantee the integrability of
$A$.  First we prove that the two conditions together imply that
$\F(A)$ is a simple foliation:
\proclaim{Lemma}
\label{lemma:transverse2:section}
For each $a\in P(A)${\rm ,} there exists
$S_a\subset P(A)$ transverse to $\F(A)${\rm ,} which intersects each leaf
of $\F(A)$ in at most one point.
\endproclaim

{\it Proof}.  The proof is a sequence of reductions and careful choices,
and is divided into several steps. So let $a\in P(\A)$ and denote by $x$
the initial point of its base path.

\vglue8pt {\it Claim}  1. We may assume that $a= 0_{x}$.
\vglue8pt 
To see this, we choose a compactly supported, time-dependent, section
$\xi$ of $\A$ so that $\xi(t, \gamma(t))= a(t)$. If $\sigma_{\xi}$
is as in the definition of admissible section (see \S
2.3), then left multiplication $T: P(\A)\to P(\A)$,
$T(b)= a_{\xi}(\t(b))b$, defines a smooth injective map.   If there is a section $S_x$ around $0_x$, as in the
statement of the lemma, it then follows that $T: S_{x}\to
P(\A)$ intersects each leaf in at most one point.  Since $S_x$ has the
same dimension as the codimension of $\F(\A)$, $S := T(S_x)$ will
be a section through $T(0_x)$ with the desired properties.  Since $T(0_x)$ and $a$ are in the same leaf, this implies the
existence of a similar section $S_a$ through $a$ (use the holonomy of the foliation $\F(\A)$ along any path from
$T(0_x)$ to $a$).
\vskip 8pt

{From} now on we fix $x\in M$ and we are going to prove the lemma
for $a= 0_x$. We also fix  local coordinates around $x$, and let $\nabla$
be the canonical flat connection on the coordinate neighborhood.  
We also choose a small neighborhood $U$ of $0_x$ in $A$ so that
the exponential map $\Exp_{\nabla}: U\to P(\A)$ is defined and is
transverse to $\F(\A)$. We are going to show that it intersects each
leaf of $\F(\A)$ in at most one point, provided $U$ is chosen
small enough.

\vglue8pt {\it Claim}  2. We may choose $U$ such that for any $v\in U\cap {\frak g}_y$ ($y\in M$)
with the property that $\Exp_{\nabla}(v)$ is homotopic to $0_{y}$, we must
have $v\in Z({\frak g}_y)$.
\pagebreak
 
Given a norm $|\cdot |$ on $\A$, the set $\set{|[v, w]|: v, w\in
{\frak g}_{y}\hbox{ with } |v|=|w|=1}$, where $y\in M$ varies in a
neighborhood of $x$, is bounded. Rescaling $|\cdot |$ if necessary, we
find a neighborhood $D$ of $x$ in $M$, and a norm $|\cdot |$ on
$A_{D}= \{ v:\pi(v)\in D\}$, such that $|[v, w]|\leq |v| |w|$ for all
$v,w\in{\frak g}_y$ with $y\in D$. We now choose $U$ so that
$U\subset A_D$ and $|v|\leq \pi $ for all $v\in U$.  If $v$ is as in
the claim, it follows from Proposition \ref{hol-hom} (see also the
proof of Lemma \ref{N-versus-tilde-N}) that parallel transport $T_{v}:
{\frak g}_{y}\to {\frak g}_y$ along the constant $\A$-path $v$
is the identity. But $T_v$ is precisely the exponential of the linear
map ${\rm ad}_{v}: {\frak g}_{y}\to {\frak g}_y$.  Since $|{\rm ad}_v|\leq
|v|\leq \pi$, it follows (see e.g.\ \cite{DoLa}) that ${\rm ad}_{v}= 0$ so
$v\in Z({\frak g}_y)$.

\vglue6pt {\it Claim}  3. We may choose $U$ such that, if $v\in U\cap {\frak g}_y$ ($y\in M$)
has the property that $\Exp_{\nabla}(v)$ is homotopic to $0_y$ then
$v= 0_y$.
\vglue6pt 
Obviously this is just a restatement of the obstruction assumptions,
combined with the previous claim.

\vglue6pt {\it Claim}  4. We may choose $U$ such that, if $v\in U$ has the property that the
base path of $\Exp_{\nabla}(v)$ is closed, then $v\in {\frak g}_y$.
\vglue6pt 
To see this, we note that the equations for the geodesics in local
coordinates reduce to:
$$ 
\left\{
\begin{array}{l}
\dot{x^i}=\sum_p b^i_p(x(t))a^p,\\
\\
\dot{a^p}=0.
\end{array}
\right.
$$ 
By the period bounding lemma ([28]),
any nontrivial periodic orbit of this system with initial condition on a
open 
set $D$ has period
$$ T\ge\frac{2\pi}{M_{D}}\qquad\hbox{ where }\quad
M_{D}=\sup_{{1\le j,k\le m\atop x\in D}}
\left\Vert\sum_p \frac{\partial b^j_p}{\partial x^k}(x)a^p\right\Vert.$$ 
Hence it suffices to make sure that $U\subset A_{D}$
where $D$ is chosen small enough so that $M_{D}< 2\pi$.
\vglue6pt

Now, for any open $O\subset P(\A)$, we consider the plaques in $O$ of
$\F(\A)$, or, equivalently, the leaves of $\F(\A)|_{O}$. For $a, b\in
O$, we write $a\sim_{O} b$ if $a$ and $b$ lie in the same plaque.
{From} now on we fix $U$ satisfying all the conditions above, and we
choose an open $O$ so that $\Exp_{\nabla}: U\to P(\A)$ intersects each
plaque inside $O$ exactly in one point. This is possible since
$\Exp_{\nabla}$ is transversal to $\F(\A)$. Apart from the pair $(O,
U)$, we also choose similar pairs $(O_{i}, U_{i})$, $i= 1, 2$, such
that $O_1O_1\subset O$, $O_{2}O_{2}\subset O_{1}$ and $O_{i}^{-1}=
O_{i}$.

\vglue6pt {\it Claim}  5. It is possible to choose a neighborhood $V$ of $x$ in $U_{2}$ so that,
for all $v\in V$,
$$  0_{y}\cdot \Exp_{\nabla}(v) \sim_{O} \Exp_{\nabla}(v) .\pagebreak$$

We know that for any $v$ there is a natural homotopy between the two
elements above. This homotopy can be viewed as a smooth map $h:
I\times U \to P(A)$ with $h(0, v)= 0\cdot \Exp_{\nabla}(v)$, $h(1,
v)= \Exp_{\nabla}(v)$, $h(t, 0_{y})= 0_{y}$. Since $I$ is compact and
$O$ is open, we can find $V$ around $x$ such that $h(I\times V)\subset
O$. Obviously $V$ has the desired property.

\vglue12pt {\it Claim}  6. It is possible to choose $V$ so that, for all $v, w\in V$,
$$ 
(\Exp_{\nabla}(v)\cdot\overline{\Exp_{\nabla}(w)}) \cdot \Exp_{\nabla}(w)
\sim_{O} \Exp_{\nabla}(v).
$$ 
\vglue12pt 
This is proved exactly as the previous claim.

\vglue12pt {\it Claim}  7. $\Exp_{\nabla}:V\to P(\A)$ intersects each leaf of $\F(\A)$ in at most
one point.
\vglue12pt 
To see this, let us assume that $v, w\in V$ have $\Exp_{\nabla}(v)\sim
\Exp_{\nabla}(w)$.  Then $a_1:= \Exp_{\nabla}(v)\cdot
\overline{\Exp_{\nabla}(w)}\in O_{1}$ will be homotopic to the trivial
$\A$-path $0_y$. On the other hand, by the choice of the pair $(O_1,
U_1)$, $a_1\sim_{O_1} \Exp_{\nabla}(u)$ for a unique $u\in U_1$. Since
$\Exp_{\nabla}(u)$ is equivalent to $0_y$, its base path must be
closed, hence, by Claim $4$ above, $u\in {\frak g}_y$.  Using Claim
$3$, it follows that $u= 0$, hence $a_{1}\sim_{O_1} 0_{y}$.  Since
$O_1O_1\subset O$, this obviously implies that
$$  a_{1}\cdot \Exp_{\nabla}(w) \sim_{O} 0_{y}\cdot\Exp_{\nabla}(w) .$$ 
Since $V$ satisfies Claims $5$ and   $6$, we get
$\Exp_{\nabla}(v)\sim_{O}\Exp_{\nabla}(w)$. Hence, by the construction
of $O$, $v= w$. This also concludes the proof of the lemma.
\hfill\qed\vglue9pt

Note that the previous lemma implies that $\G(A)$ has a natural
quotient differentiable structure: the charts are just the $S_{a}$'s,
and the change of coordinates is smooth since it is just the holonomy
of $\F(A)$. Hence we can complete the proof of Theorem \ref{main thm}
by showing that:

\proclaim{Lemma}
\label{lemma:Lie:groupoid}
For the quotient differentiable structure{\rm ,} $\G(A)$ is a Lie groupoid
with Lie algebroid $A$.
\endproclaim 

\demo{Proof} 
It is clear from the definitions that $A$ can be identified with
$T^{\s}_M \G(A)$ and that under this identification $\#$ coincides
with the differential of the target~$\t$. So we need only to check
that the bracket of right-invariant vector fields on $\G(A)$ is
identified with the bracket of sections of $A$. For this we note that on
one hand, the bracket is completely determined by the infinitesimal
flow of sections through the basic formula (\ref{Lie-flows}). On the
other hand, we now know that the exponential
$\exp:\Gamma(A))\to\Gamma(\G(A))$ is injective in a neighborhood of
the zero section, and so Proposition \ref{prop:exponential} shows that
the infinitesimal flow of a section $\al$ is the infinitesimal flow of
the right-invariant vector field on $\G(A)$ determined by
$\al$. Hence, we must have $\A(\G(A))=A$.
\enddemo

{\it Remark} 4.12.
\label{one:in:a:leaf}
The proof above (namely an argument similar to Claim $1$ above) shows
that, in the main theorem, it suffices to require that for each leaf~$L$,
there exists $x\in L$ satisfying the two obstructions.

\section{Examples and applications}  %
\label{Examples/Applications}        %

In this section we review the known integrability criteria,
we derive them from Theorem \ref{main thm}, and present an application
to the theory of transversally parallelizable foliations.

\vglue6pt
 5.1. {\it  Local integrability}.
Regarding the local nature of integrability, with respect to the base
manifold $M$, note that
\begin{itemize}
\item From Examples \ref{ex:second:obstruction} and
\ref{example:Weinstein} we learn that a Lie algebroid can be locally
integrable (i.e.\ each point has a neighborhood $U$ so that $A|_{U}$ is
integrable), and not globally integrable. This shows that the
integrability problem is not a local one.
\item Example \ref{example:second:obstruction:2} shows that there are
algebroids which are not even locally integrable.
\end{itemize}
However, a general ``local integrability'' result has long been
assumed to be true, namely the integrability by {\it local
groupoids}. This result was first announced by Pradines, but a proof
has never been published.  One of the main difficulties is that, if
one tries to extend the known result from Lie groups (see
e.g. \cite{DuKo}), one faces the problem of finding a
CBH-formula. However, with the Weinstein groupoid at hand (and its
description as a leaf space) this result can be proved quite
easily. 

For a {\it local Lie groupoid} the structure maps are only defined on
(and the usual properties only hold for) elements which are close
enough to the space $M$ of units (these are obvious generalizations of
Cartan's local Lie groups, as explained in Section $1.8$ of
\cite{DuKo}).

\proclaim{{C}orollary} 
Any Lie algebroid is integrable by a local Lie groupoid.
\endproclaim 
\vglue-28pt
\phantom{bugger}

{\it Proof}.  One uses exactly the same arguments as in the proof
of Claim $2$ and Claims $4$ through $6$ in Lemma
\ref{lemma:transverse2:section}, namely: We choose a connection
$\nabla$ on~$\A$, a neighborhood $U$ of $M$ in $A$ and an open $O$ in
$P(\A)$, with $O= O^{-1}$, such that $\Exp_{\nabla}:U\to P(\A)$
intersects each plaque of $\F(\A)$ in $O$ in exactly one point.
Eventually choosing smaller pairs $(O, U)$ (similar to the $(O_i,
U_i)$ in the cited proof), the structure of local groupoid will be
defined on $U$: the inverse of $v\in U$ is the unique $\bar{v}\in U$
with the property that $\overline{\Exp_{\nabla}(v)}\sim_{O}
\Exp_{\nabla}(\bar{v})$; the multiplication $v\cdot w$ of $v, w\in U$
is defined only for pairs $(u, v)$ for which
$\Exp_{\nabla}(v)\Exp_{\nabla}(w)\in O$, and is the unique element with
the property that the last product of exponentials is $\sim_{O}
\Exp_{\nabla}(v\cdot w)$. The associativity around the units is proved
exactly as Claims $5$ and $6$ of the cited lemma, while the fact
that the resulting local groupoid integrates $A$ is a variation of
Lemma \ref{lemma:Lie:groupoid}. \phantom{whereis}
\hfill\qed
\vglue9pt
 
5.2. {\it Integrability criteria}.  We start with the following general integrability criterion, which is an
obvious consequence of our main result, and which implies most of the known
results (or even much stronger versions of them):

\proclaim{{C}orollary}
\label{general:criteria}
If $N_{x}(A)$ is trivial for all $x\in M${\rm ,} then $A$ is integrable. In
particular{\rm ,} $\A$ is integrable if any of the following three
conditions holds for all $x\in M${\rm :}
\begin{itemize}
\item[{\rm (i)}] the Lie algebras $\gg_x$ are semi\/{\rm -}\/simple {\rm (}\/more generally{\rm ,} if
they have trivial center\/{\rm );}
\item[{\rm (ii)}] the leaves $L_x$ are {\rm 2-}\/connected {\rm (}\/more generally{\rm ,} if
$\pi_{2}(L_{x})$ have only elements of finite order\/{\rm );}
\item[{\rm (iii)}] there is a splitting $\sigma: TL_x\to A|_{L_x}$ of the
anchor{\rm ,} which is compatible with the Lie bracket\/{\rm ;}
\end{itemize}

\endproclaim 

We now briefly deduce the known integrability results.

\vglue12pt 5.2.1. {\it Lie algebra bundles}. For Lie algebroids with zero
anchor map, the orbits are the points of $M$, so the conditions of the
main theorem are trivially satisfied and we obtain the following
result of Douady and Lazard \cite{DoLa}:

\proclaim{{C}orollary}
Any Lie algebra bundle is integrable to a Lie group bundle.
\endproclaim 

5.2.2. {\it Transitive algebroids}.
In the case of Lie algebroids with surjective anchor the main theorem
becomes
(see also Remark \ref{one:in:a:leaf}):

\proclaim{{C}orollary} Let $\A$ be a transitive Lie algebroid over $M$. Then $\A$ is
integrable if and only if $N_{x}(A)$ is discrete in $A_x$ for one
{\rm (}\/or{\rm ,} equivalently{\rm ,} all\/{\rm )} $x\in M$.
\endproclaim 

We mention in passing that this is strongly related to Mackenzie's
criteria \cite{Mack1}, and the interested reader will be able to provide the
precise relation with his cohomological obstruction.

There are some obvious consequences of this result. For example:

\proclaim{{C}orollary}
Every transitive Lie algebroid $\A$ over a $2$\/{\rm -}\/connected base $M$ is
integrable. 
\endproclaim 

Note also that since $\s^{-1}(x)$ is a principal $\G(\A)_x$-bundle
over $M$, it follows that if $M$ is contractible then $A$ is in fact
isomorphic to a direct sum $TM\oplus \gg$ (compatible with the Lie
brackets), where $\gg= {\frak g}_x$.  Hence:

\proclaim{{C}orollary} Any transitive Lie algebroid over a contractible base $M$
is isomorphic to $TM\oplus \gg$ for some Lie algebra ${\frak g}$.
\endproclaim 

In Mackenzie's approach this result is first obtained in order to
to construct his obstruction.

\vglue12pt 5.2.3. {\it Regular Lie algebroids}.
Although many of the known integrability criteria require regular
algebroids, it turns out that regularity is superfluous (see below).
This is the case, for example, with the Dazord-Hector (\cite{DaHe})
integrability criteria for totally aspherical regular Poisson
manifolds, and with Nistor's results \cite{Ni} on the integrability of
regular algebroids with anchor possessing a splitting compatible with
the Lie bracket, or with semi-simple kernels.

Let us mentioned, however, a result which fails in the nonregular
case as shown by Example \ref{example:second:obstruction:2}:

\proclaim{{C}orollary} 
Any regular Lie algebroid is locally integrable.
\endproclaim 

This follows because regular foliations are locally trivial. As in
the transitive case, it is possible to describe explicitly
the local structure of regular algebroids. Choosing local coordinates
in $M$ so that the foliation becomes the obvious $p$-dimensional
foliation on ${\Bbb R}^p\times {\Bbb R}^q$, we then obtain that,
locally,
the algebroid is $T{\Bbb R}^p\times {\frak g}$, where
${\frak g}$ is a bundle of Lie algebras over ${\Bbb R}^q$.

\vglue12pt 5.2.4. {\it Semi\/{\rm -}\/direct products}.
Closely related to Palais' integrability \cite{Palais} of infinitesimal
actions of Lie algebras ${\frak g}$ on manifolds $M$
is the integrability of the transformation Lie algebroid
$\A= {\frak g}\times M$. Recall that, as a vector bundle,
$A$ is just the trivial vector bundle with fiber ${\frak g}$,
the anchor is the infinitesimal action, while the bracket on
$\Gamma(A)=C^{\infty}(M; {\frak g})$ is uniquely determined by the
Leibniz rule and the Lie bracket of ${\frak g}$.

If $G_x\subset {\cal G}(\frak g)$ is the connected Lie group with Lie algebra ${\frak g}_x$, $N_x(\A)$ sits inside
$\pi_1(G_x)$, so the 
conditions of the main theorem are satisfied; hence

\proclaim{{C}orollary}
For any infinitesimal action of the Lie algebra ${\frak g}$ on $M${\rm ,} the
transformation Lie algebroid ${\frak g}\times M$ is integrable.
\endproclaim

This is known as Dazord's criterion (cf.~\cite{CaWe}), but it also
appears implicitly in Palais' work \cite{Palais}. Implicit in Palais'
work is also the precise relation between this result and the
integrability of infinitesimal actions. This relation has been clearly
explained by Moerdijk-Mr\v{c}un in \cite{MoMr}, where the reader can
find various extensions \pagebreak to semi-direct products of algebroids.  Let us
point out that exactly the same argument as above shows that the
semi-direct product of an integrable algebroid by a regular foliation
is integrable, and this is one of the main results of \cite{MoMr}.


\vglue5pt 5.2.5. {\it Algebras of vector fields and quasi\/{\rm -}\/foliations}.
For a Lie algebroid $A$ over $M$ we say that the the anchor is
{\it almost injective} at $x_0\in M$ if there is a neigborhood $U$ of
$x_0$ in $M$ and an open dense subset $O$ of $U$ such that $\#_x$ is
injective for all $x\in O$. Note that if the anchor is injective at
$x_0$ then it is almost injective at $x_0$. We say that the anchor is
almost injective if it is almost injective at every point.

Any Lie subalgebra $\Gamma$ of $\X(M)$ which is a finitely generated
projective $C^{\infty}(M)$-module is the space of sections of an
algebroid with almost injective anchor.  This
produces a large class of examples of Lie algebroids, including all
regular foliations.  As explained in \cite{Ni}, such $\Gamma$'s arise
naturally in the analysis on manifolds with corners as algebras of
vector fields with a certain behavior on the faces of $M$. Their
integrability is relevant to various aspects of analysis and
quantization (see \cite{NiXuWe} for details).  Such algebroids were
also studied by Claire Debord in her Ph.D.~thesis (\cite{Deb}), and
they give rise to {\it quasi-foliations} of $M$.

Our main result implies the following integrability criterion due to Debord:
\vglue-16pt
\phantom{du}

\proclaim{{C}orollary}
\label{cor:almost:injective}
A Lie algebroid with injective anchor on a dense open set is
integrable.
\endproclaim 
\vglue-16pt
\phantom{du}

 This result follows from the following general continuity type
property of the monodromy groups, combined with our main
result.  

\proclaim{Proposition}
\label{prop:1}
Let $A$ be a Lie algebroid over $M$. Then{\rm ,} for every $x\in M${\rm ,} there
exists an open set $U\subset A$ containing $0_x$ such that $N(A)\cap
U$ has enough smooth local sections{\rm ,} i.e.{\rm ,} for all $v\in N(A)\cap U${\rm ,}
there is a section $\al:V\to A$ defined in a neighborhood $V$ of
$\pi(v)\in M$ such that $\al(\pi(v))= v$ and $\alpha(y)\in
N_{y}(A)\cap U$ for all $y\in V$.
\endproclaim

{\it Proof}.
Using Claim $2$ and Claim $4$ in the proof of the main theorem (which are 
independent of any integrability conditions!), we find a neighborhood $U$ of 
$0_x$ in $A$ such that:
\begin{itemize}
\item[(a)] $\Exp_{\nabla}: U\to P(A)$ is well defined and transversal
  to the foliation;
\item[(b)] if $w\in U$ has $\Exp_{\nabla}(w)\sim 0$, then $w\in Z(\gg)$.
\end{itemize}
Here, as in the proof of the main theorem, $\nabla$ is a fixed local
connection around~$x$.
 
We claim that $U$ does satisfy the desired property. Assume $v\in
N_{z}(A)\cap U$. Then a homotopy between the zero path and $v$ can be
viewed as a leafwise path in $P(A)$, connecting the constant path
$0_z$ with the constant path $v$. We denote by $H$ the induced
holonomy, where we use $(\Exp_{\nabla}, 0_z)$ as transversal through
$0_z$, and $(\Exp_{\nabla}, v)$ as transversal through $v$. Hence,
$H:(U_0, 0_z)\to (U_1, v)$ for some neighborhoods $U_0$ and $U_1$ in
$U$, and $\Exp_{\nabla}(H(w))\sim \Exp_{\nabla}(w)$ for all $w\in
U_0$. From condition (b) it follows that $H(0_y)$ is a constant
$A$-path in $Z(\gg_y)$ for all $y$ in $V=\{y:0_y\in U_0\}$. Therefore we
can choose $\al: V\to A$ to be the section defined by $\al(y)\equiv H(0_y)$.
\hfill\qed\vglue12pt

Let us point out also the following improvement of Corollary 5.9:

\proclaim{{C}orollary}
A Lie algebroid $A$ such that the monodromy groups
$N_{x}(A)$ are trivial for $x$ on a dense open subset of the
base manifold $M$ is integrable.
\endproclaim

5.2.6. {\it Poisson manifolds}.
The Weinstein groupoid of the algebroid associated to a Poisson
manifold (the cotangent bundle $T^*M$) is precisely the phase space
$\G$ of the Poisson sigma-model studied by Cattaneo and Felder in
\cite{CaFe}. Our constructions explain the constructions in
\cite{CaFe}, while our main result clarifies the smoothness of the
Poisson-sigma model $\G$.

The following obvious application of our general criteria, is the main
positive result of \cite{CaFe}.

\proclaim{{C}orollary}
Any Poisson structure on a domain in ${\Bbb R}^2$ is integrable.
\endproclaim 
\phantom{huh}
\vglue-18pt 

The result is certainly not true in higher dimension, as shown by\break
Weinstein's example of a nonintegrable regular Poisson structure in
$\Rr^3-0$ (Example~\ref{example:Weinstein}).

In general, our main result applied to this context describes the
precise obstructions for the integrability of Poisson manifolds. Let
us point out the following simple integrability result:

\proclaim{{C}orollary}
Any Poisson manifold where the symplectic leaves have vanishing second
homotopy groups is integrable.
\endproclaim 

The integrability criterion of Dazord and Hector (\cite{DaHe}) is in fact
this result specialized to the case of a regular Poisson manifold for
which the foliation has no vanishing cycles.

Note also that the monodromy groups of the regular symplectic leaves
$L$ (i.e.\ around which the rank is locally maximal) of a Poisson
manifold $M$ are particularly simple, as it is the associated
monodromy map
$$  \partial: \pi_{2}(L, x)\to \NN^{*}_{x}(L) .$$ 
Indeed, since the kernel $\NN^{*}(L)$ of $\#$ over $L$ is abelian, by
Lemma \ref{compute} we can use any linear splitting $\sigma$. The
resulting cohomology class
$$  \Omega_{L}= [\Omega_{\sigma}]\in H^{2}(L; \NN^{*}(L)) $$ 
is independent of the splitting $\sigma$, $\partial$ is just the
integration of $\Omega_{L}$ over elements in $\pi_{2}(L, x)$, and its
image defines the monodromy groups
$$  N_{x}\subset \NN^{*}_{x}(L) .$$ 
These groups appeared already in the work of Alcalde-Cuesta and Hector
(\cite{AlHe}), who also recognized the two obstructions to
integrability (in an equivalent form). To our knowledge, their work
contained the most detailed investigation up to this date of the
integrability problem for Poisson manifolds. Unfortunately, their
methods do not seem to apply beyond the regular case.

Notice also that if $M$ is regular and $\F$ is its symplectic foliation,
then using a global splitting $\sigma$ for $\#$ one gets a globally
defined cohomology class $\Omega\in H^{2}(\F; \NN^{*})$ which lies in
the foliated cohomology with coefficients in the kernel of
$\#$. Clearly, $\Omega_{L}= \Omega|_{L}$ for each $L$. For more on
integrability of Poisson manifolds we refer the reader for the
upcoming article \cite{CraFer}.

\vglue12pt 5.2.7. {\it Van Est\/{\rm '}\/s argument}.
Probably the most elegant proof of the integrability
of Lie algebras is Van Est's cohomological argument
which we briefly recall. Given a Lie algebra ${\frak g}$,
we form the exact sequence
$0\to Z({\frak g})\to {\frak g}\to {\rm ad}({\frak g})\to 0$.
Here ${\rm ad}({\frak g})$ is easily seen to be integrable
(it is a Lie sub-algebra of ${\frak g}{\frak l}({\frak g})$).
Also recall that simply connected Lie groups are automatically\break
$2$-connected.
The core of Van Est's argument is then the following result
for the particular case of Lie algebras

\proclaim{{C}orollary} If $B$ fits into an exact sequence of Lie algebroids
$$  0\to E\to B \stackrel{\pi}{\to}A \to 0 $$ 
with $E$ abelian{\rm ,} and $A$ integrable by a groupoid with $2$\/{\rm -}\/connected
$\s$\/{\rm -}\/fibers{\rm ,} then $B$ is integrable.
\endproclaim

This result for Lie algebroids is Theorem $5$ of
\cite{Cra}. Interestingly enough, it shows that the integrability
criterion of Dazord and Hector \cite{DaHe} mentioned above is actually
Van Est's argument applied to regular Poisson manifolds.

The proof in \cite{Cra} is an extension of Van Est's cohomological methods.
Using a splitting $\sigma$ of $\pi$ we obtain an action of $A$
on $E$, and a $2$-cocycle $\Omega_{\sigma}$ on $A$ with values
on $E$. This is well known (see e.g.\ \cite{Mack1}), and can also be viewed
as 
an extension of the constructions in Section 3.3.
We can then form the group of periods ${\rm Per}_{x}\subset E_x$ of
$\Omega_{\sigma}$.
The cohomological proof actually shows that $B$ is integrable
provided $A$ is, and the groups $P_x$ vanish (cf.~Remark $5$ and
Corollary~$2$ in \cite{Cra}).

Let us briefly point out how our result implies (and further clarifies)
the previous corollary.
Let $x\in M$ sitting in a singular leaf $L$.
The necessary information is organized in the following diagram

\figin{41}{1000}%

\noindent Here $\s_{A}$, $\t_A$, ${\frak g}(A)$, $\partial_{A}$ are respectively
the source and target map, the kernel of the anchor, and the monodromy map
of $A$, and we use analogous notations for $B$. Also,
$\partial$ is the boundary
map in homotopy associated to $\s_{A}^{-1}(x)\to L$ with
fiber $\G(\A)_{x}$ and $j$ is the obvious inclusion for which the
image is precisely $\tilde{N}_{x}(A)$. Finally,
$\partial_{E}$ denotes the monodromy map associated to the
exact sequence in the corollary (constructed exactly as the monodromy
map of \S 3.2), and its image
is precisely the group of periods ${\rm Per}_x$.

\proclaim{Lemma}
There is a short exact sequence of abelian groups\/{\rm :}\/
$$  0\to {\rm Per}_{x} \to \tilde{N}_{x}(B) \to \tilde{N}_{x}(A)\to 0 .$$ 
\endproclaim 
\vglue-20pt
{\it Proof}.
These follows by diagram chasing since the two horizontal sequences
above are exact.
\hfill\qed\vglue18pt

Therefore $\tilde{N}_{x}(B)$ appears as a twisted semi-direct product
of $\tilde{N}_{x}(A)$ and ${\rm Per}_x$.  The simplest case where our main
theorem applies is when ${\rm Per}_{x}$ vanishes. This gives precisely the
corollary (and its stronger version) above.

\demo{{\rm 5.3.} Transversally parallelizable foliations} 
Historically, the first examples of nonintegrable Lie algebroids
\cite{Mol} came from Molino's treatment (see \cite{Molino}) of
trans\-ver\-sa\-lly parallelizable foliations which we now briefly
recall.

Given a foliation $\F$ of $M$, let us denote by $l(M, \F)$ the algebra
of transversal vector fields, i.e.\ sections of the normal bundle
which can be locally projected along submersions which locally define
the foliation. Then $(M, \F)$ is {\it transversally parallelizable} if
its normal bundle admits a global frame consisting of transversal
vector fields.  In this case the Lie algebra $l(M, \F)$ is free as a
module over the space $\Omega_{b}^{0}(M, \F)$ of basic functions, on
which it acts by derivations.

Let us show that the Lie bracket on $l(M, \F)$ is of the type studied
in this paper. We assume for simplicity that $M$ is compact. Then the
closures of the leaves of $\F$ form a new foliation $\bar{\F}$ on $M$,
with leaf space a smooth (Hausdorff) manifold $W=M/\bar{\F}$, called
the {\it basic manifold} of the foliation.  Since $\F$ and $\bar{\F}$
have the same basic functions, $l(M, \F)$ is the space of sections of
a transitive Lie algebroid over $W$, which we denote by $A(M,
\F)$. Its anchor $\#$ is just the action of $l(M, \F)$ on
$\Omega_{b}^{0}(M, \F)\cong C^{\infty}(W)$, and the kernel of $\#$ has
the following geometric interpretation.  For each leaf $L$ of $\F$,
the foliation $(\bar{L}, \F|_{\bar{L}})$ is transversally
parallelizable with dense leaves. It follows that $l(\bar{L},
\F|_{\bar{L}})$ is a finite-dimensional Lie algebra, and moreover,
$(\bar{L}, \F|_{\bar{L}})$ is a Lie foliation induced by a canonical
$l(\bar{L}, \F|_{\bar{L}})$-valued Maurer-Cartan form. Denoting by
$w\in W$ the point defined by $\bar{L}$, we   see that $l(\bar{L},
\F|_{\bar{L}})$ is canonically isomorphic to ${\rm Ker}(\#_{w})$.  This
shows that all the Lie algebras $l(\bar{L}, \F|_{\bar{L}})$ are
isomorphic. The resulting Lie algebra ${\frak g}(M, \F)$ (defined
up to isomorphisms) is usually called the {\it structural Lie
algebra} of the foliation.

The main result of Almeida and Molino in \cite{Mol} says that $(M, \F)$
is developable (i.e.\ its lift to the universal cover of $M$ is
simple) if and only if the Lie algebroid $A(M, \F)$ is integrable.
This discussion extends to transversally complete foliations $(M, \F)$
without any compactness assumption on $M$ (see \cite{Molino}).

Now, our constructions produce a monodromy map $\partial : \pi_{2}(W)
\to G(M, \F)$ with values in the simply connected Lie group
integrating the structural Lie algebra ${\frak g}(M, \F)$, which
controls the developability of the foliation:

\proclaim{{C}orollary}
A transversally parallelizable foliation $(M, \F)$ on a compact
manifold $M$ is developable if and only if the image of the monodromy
map
$$  \partial : \pi_{2}(W) \to G(M, \F) $$ 
is discrete. 
\endproclaim 

A simple consequence of this result is:

\proclaim{{C}orollary} Let $(M, \F)$ be a transversally parallelizable manifold
on a compact manifold  $M$. Then $(M, \F)$ is developable in any of the following cases\/{\rm :}
\begin{itemize}
\item[{\rm (i)}] the structural Lie algebra ${\frak g}(M, \F)$
has trivial center\/{\rm ;}\/
\item[{\rm (ii)}] $\pi_{2}(W)$ has only elements of finite order.
\end{itemize}
\endproclaim 

This result should be compared with Corollary 1,  p.~301, and Corollary
1  p.~303 in \cite{Molino}.

 \vglue16pt \centerline{\bf Appendix A.  Flows}
\vglue12pt

In this appendix we discuss the flows associated to sections of Lie
algebroids, which generalize the ordinary flows of vector fields (sections
of
$A=TM$). This is used throughout the paper, most notably for defining the
equivalence relation on $\A$-paths (\S 1.3).
As in the main body of the paper, $\A$ denotes a Lie algebroid over $M$,
$\#: \A\to TM$ denotes its anchor and $\pi: \A\to M$ the projection.

\demo{{\rm A.1.} Flows and infinitesimal flows} Given a time-dependent vector field $X$ on $M$, we denote by
$\Phi_{X}^{t, s}$ its flow from time $s$ to time $t$. Hence
$$  \frac{d}{dt} \Phi_{X}^{t, s}(x)= X(t, \Phi_{X}^{t, s}(x)),
\ \ \Phi_{X}^{s, s}(x)= x  .$$ 
We have $\Phi^{t, s}\Phi^{s, u}= \Phi^{t, u}$ and, when $X$ is
autonomous, $\Phi_{X}^{t, s}= \Phi_{X}^{t- s}$ only depends on $t-s$.
Differentiating, we obtain the infinitesimal flow of $X$:
$$ 
\phi_{X}^{t, s}(x)\equiv (d\Phi_{X}^{t, s})_{x}: T_x M\to
T_{\Phi_{X}^{t, s}(x)}M.
$$ 

Let us assume now that $\G$ is a Lie groupoid integrating the
algebroid $\A$. Given a time-dependent section $\alpha$ of $\A$, we
denote by the same letter the right invariant (time-dependent)
vector field on $\G$ induced by $\alpha$, and by
$\varphi_{\alpha}^{t, s}: \G\to \G$ its flow. If $x=\s(g)$ and
$y=\t(g)$, then $\varphi_{\alpha}^{t, s}(g)$ is the arrow
$$  \varphi_{\alpha}^{t, s}(g): x\rmap \Phi_{\#\al}^{t, s}(y)$$ 
and also satisfies the right-invariance property:
$$  \varphi_{\alpha}^{t, s}(g)= \varphi_{\alpha}^{t, s}(y)g.$$ 
The {\it infinitesimal flow} of $\alpha$,
$$  \phi_{\alpha}^{t, s}(x): A_x\to   A_{\Phi_{\#\al}^{t, s}(y)},$$ 
is defined as
\begin{equation}
\label{flow-G-A}
 \phi_{\alpha}^{t, s}(x)\equiv
          (d R_{\varphi_{\alpha}^{s, t}(x)})_{\varphi_{\alpha}^{t, s}(x)}
          (d\varphi_{\alpha}^{t, s})_{x}. \speqnu{A.1}
\end{equation}
The classical relation between Lie brackets and flows translates at this
level
to
\begin{equation}
\label{Lie-flows}
{\left. \frac{d}{dt}\right|}_{t= s} (\phi_{\alpha}^{t, s})^{*}\beta =
[\alpha^s, \beta], \speqnu{A.2}
\end{equation}
where we have set
$$  (\phi_{\alpha}^{t, s})^{*}(\beta)(x)=
   \phi_{\alpha}^{s, t}\beta(\phi_{\#\alpha}^{t, s}(x)).$$ 

We wish to extend the infinitesimal flow to sections of general Lie
algebroids, not necessarily integrable. For this we can use the
following general construction of (infinitesimal) flows. Let us assume
that $E$ is a vector bundle over $M$.  A derivation on $E$ is a pair
$(D,X)$ where $D: \Gamma(E)\to \Gamma(E)$ is a differential operator,
$X$ is a vector field on $M$, satisfying the Leibniz rule
$$  D(f\alpha)= fD(\alpha)+ X(f)\alpha, \ \
     \hbox{for all } \ f\in C^{\infty}(M), \alpha\in \Gamma(E).$$ 
Now, any time-dependent derivation $(D, X)$ on $E$ has an associated
(infinitesimal) flow: a standard argument shows that there is a family
of linear isomorphisms
$$  \phi_{D}^{t, s}(x): E_{x}\to E_{\Phi_{X}^{t, s}(x)} , \pagebreak$$ 
which is characterized uniquely by the properties
\begin{itemize}
\item[(a)] $\phi^{t,s}_{D}\phi^{s, u}_{D}= \phi^{t, u}_{D}$,
$\phi^{t, t}_{D}= {\rm Id}$;
\item[(b)]
${\left. \frac{d}{dt}\right|}_{t= s} (\phi_{D}^{t, s})^{*}\beta
=D^{s}(\beta)$,
for all sections $\beta\in \Gamma(E)$.
\end{itemize}
Here $D^{t}$ is $D$ at the fixed time $t$, and
$(\phi_{D}^{t,s})^{*}\beta = \phi_{D}^{s, t}\beta\Phi_{X}^{t, s}$.

Alternatively, one can use the groupoid ${\rm Aut}(E)$ over $E$, for which the
arrows from $x$ to $y$ are all linear isomorphisms $E_{x}\to E_{y}$.
Its Lie algebroid is usually denoted by $DO(E)$, and its sections are
precisely derivations of $E$ (cf.~\cite{HiMa}, \cite{Mack1}). Hence $(D,X)$
can be viewed as a time-dependent section of $DO(E)$, and then
$\phi_{D}^{t, s}$ is just the associated flow on ${\rm Aut}(E)$. Both
definitions of $\phi_{D}^{t, s}(x)$ show they are defined whenever
$\Phi_{X}^{t, s}(x)$ is defined.

Most flows in differential geometry (e.g.\ the flows of vector fields,
parallel transport, etc.) are obtained in this way.
\enddemo

A.2. {\it The infinitesimal flow of a section}. We apply the previous construction to a time-dependent section
$\alpha$ of the Lie algebroid $\A$, where $X=\#\al$ and $D=
[\alpha, -]: \Gamma(\A)\to \Gamma(\A)$. The resulting flow
$$  \phi_{\alpha}^{t, s}(x): A_x\to   A_{\Phi_{\#\al}^{t, s}(x)}$$ 
is uniquely determined by $\phi^{t,s}_{\alpha}\phi^{s, u}_{\alpha}=
\phi^{t, u}_{\alpha}$, $\phi^{t, t}_{\alpha}= {\rm Id}$, and the formula
(\ref{Lie-flows}) above.  In particular, if $\A$ is integrable, then
$\phi_{\alpha}^{t, s}$ coincides with (\ref{flow-G-A}) above. As in
the case of vector fields, if $\alpha$ is autonomous, then
$\phi_{\alpha}^{t, s}= \phi_{\alpha}^{t- s}$ only depends on $t-s$.

Let us indicate an alternative description. Recall that
on $A^*$ one has a Poisson bracket $\set{~,~}_A$ which is
linear on the fibers. A section $\al$ of $A$ defines in a natural
way a function $f_\al: A^*\to \Rr$ which is linear on the fibers
(``evaluation''), and we denote by $X_{\al}$ the Hamiltonian vector
field associated with $f_\al$. It is easy to see
(cf.~\cite{Fer1}) that:
\begin{itemize}
\item[(a)] The assignment $\al\mapsto f_\al$ defines a Lie algebra
homomorphism $(\Gamma(A),[~,~])\break\to (C^\infty(A^*),\set{~,~}_A)$;
\item[(b)] $X_{\al}$ is $\pi$-related to $\#\al$:
$\pi_*X_{\al}=\#\al$, where $\pi: A^*\to M$ is the natural
projection;
\end{itemize}
For each $t$, the flow $\Phi^{s,t}$ of $X_{\al}$ defines a Poisson
automorphism of $A^*$ (wherever defined), which maps linearly fibers
to fibers of $A^*$.  So, in fact, $\Phi^{s,t}:A^*\to A^*$ is a bundle map
and from (b) we have that it covers $\Phi^{s,t}_{\#\al}$, the flow of
$\#\al$. By transposition we obtain the infinitesimal flow
$\phi_{\alpha}^{t, s}(x): A_x\to   \pagebreak A_{\Phi_{\#\al}^{t, s}(x)}$.

\demo{{E}xample {\rm A.1}} As a simple example, consider a Lie algebra $\A= \gg$ as a Lie
algebroid over a point, and $\alpha\in \gg$ (a constant section). The
Poisson structure on the dual $\gg^*$ is the Kirillov-Kostant bracket
and so the hamiltonian flow on $\gg^*$ of the evaluation function
$f_\al$ is given by the co-adjoint action. It follows that the
infinitesimal flow of $\al$ is then $\phi_{\alpha}^{t}=
{\rm Ad}(\exp(t\alpha))$.
\enddemo

This example shows that one can think of the infinitesimal flow of a
section as a generalization of the adjoint action, although for a
general Lie algebroid it does not make sense to speak of the adjoint
representation!


\begin{references} 
 
\bibitem{AlHe} \name{F.~Alcade-Cuesta} and \name{G.~Hector}, Int\'egration symplectique
des vari\'et\'es de Poisson r\'eguli\`eres, {\it Israel
J.~Math.}~{\bf 90}
(1995), 125--165.

\bibitem{Mol} \name{R.~Almeida} and \name{P.~Molino}, Suites d'Atiyah et
feuilletages transversalement complets, {\it C.~R.~Acad.~Sci.~Paris
S\'er.~I Math.~}{\bf 300} (1985), 13--15.

\bibitem{CaWe} \name{A.~Cannas da Silva} and \name{A.~Weinstein}, {\it Geometric
Models for Noncommutative Algebras}, {\it Berkeley Math.\  Lectures\/},
{\bf 10}, A.\ M.\ S., Providence, RI,  1999.

\bibitem{CaFe} \name{A.\ S.~Cattaneo} and \name{G.~Felder}, Poisson sigma models and
  symplectic groupoids, in {\it Quantization of Singular Symplectic
  Quotients}  (N.~P.~Landsman, M.~Pflaum, and M.~Schlichenmeier, eds.),
{\it Progr.\ in Math\/}.\ {\bf 198} (2001), 41--73.

\bibitem{Cra} \name{M.~Crainic}, Differentiable and algebroid cohomology, Van
Est isomorphisms, and characteristic classes, preprint
  {\it math.DG/0008064},
{\it Comment.\ Math.\ Helv\/}., to appear.

\bibitem{Cra2} \bibline, Connections up to homotopy and
characteristic classes, preprint {math.DG/ 0010085}.

\bibitem{CraFer} \name{M.~Crainic} and \name{R.~L.~Fernandes}, Integrability of
Poisson manifolds, preprint  math.DG/ 0210152.

\bibitem{DaHe} \name{P.~ Dazord} and \name{G.~Hector}, Int\'egration symplectique des
vari\'et\'es de Poisson totalement asph\'eriques, in {\it Symplectic
Geometry, Groupoids and Integrable Systems} (Berkeley, CA, 1989),  {\it MSRI
Publ\/}.\ {\bf 20}
(1991), 37--72.

\bibitem{Deb} \name{C.~Debord}, Feuilletages singulaires et
groupo\"{\i}des d'holonomie, Ph.~D.~Thesis, Universit\'e Paul
Sabatier, Toulouse, 2000.

\bibitem{DoLa} \name{A.~Douady} and \name{M.~Lazard}, Espaces fibr\'es en alg\`ebres
de Lie et en groupes, {\it Invent.~Math.}~{\bf 1} (1966), 133--151.

\bibitem{DuKo} \name{J.~Duistermaat} and \name{J.~Kolk}, {\it Lie Groups},
{\it Universitext\/}, Springer-Verlag, New York,   2000.

\bibitem{Fer2} \name{R.~L.~Fernandes}, Connections in Poisson geometry
I.\ Holonomy and invariants, {\it J.~Differential Geom.}\ {\bf 54}
(2000), 303--365.


\bibitem{Fer1}  \bibline, Lie algebroids, holonomy and
characteristic classe,  
{\it Adv.~in Math.\/} {\bf 170} (2002), 119--179.


\bibitem{Ginz} \name{V.~Ginzburg}, Grothendieck groups of Poisson vector
bundles, {\it J.\ Symplectic Geometry\/} {\bf 1} (2001), 121--169.

\bibitem{HiMa} \name{P.\ J.~Higgins} and \name{K.~Mackenzie}, Algebraic constructions
in the category of Lie algebroids, {\it J.~Algebra}~{\bf 129}
(1990), 194--230.

\bibitem{Lang} \name{S.~Lang}, {\it Differential Manifolds},
Second edition, Springer-Verlag, New-York, 1985.

\bibitem{Mack1} \name{K.~Mackenzie}, {\it Lie Groupoids and Lie Algebroids
in Differential Geometry}, {\it London Math.\ Soc.~Lecture Notes
Ser\/}.\
{\bf 124}, Cambridge Univ.~Press, Cambridge, 1987.

\bibitem{Mack2} \bibline, Lie algebroids and Lie pseudoalgebras,
{\it Bull.~London Math.~Soc.~}{\bf 27} (1995), 97--147.

\bibitem{MaXu} \name{K.~Mackenzie} and \name{P.~Xu}, Integration of Lie
bialgebroids, {\it Topology} {\bf 39} (2000), 445--467.


\bibitem{MoMr} \name{I.~Moerdijk} and \name{J.~Mr\v{c}un}, On integrability of
infinitesimal actions, {\it Amer.\ J.\ Math\/}.\ {\bf 124} (2002),
567--593.




\bibitem{Molino} \name{P.~Molino}, \'Etude des feuilletages transversalement
complets et applications, {\it Ann.~Sci.\ \'Ecole Norm.~Sup.}~{\bf 10}
(1977), 289--307.

\bibitem{Ni} \name{V.~Nistor}, Groupoids and the integration of Lie algebroids,
{\it J.~Math.~Soc.~Japan} {\bf 52} (2000), 847--868.

\bibitem{NiXuWe} \name{V.~Nistor, A.~Weinstein}, and \name{P.~Xu}, Pseudodifferential
operators on differential groupoids,
{\it Pacific~J.~Math.}~{\bf 189} (1999), 117--152.

\bibitem{Palais} \name{R.~Palais}, A global formulation of the Lie theory of
transformation groups, {\it Mem.~Amer.\ Math.\ Soc.}~{\bf 22} (1957),
Providence, RI.

\bibitem{Sing} \name{I.\ M.~~Singer} and \name{S.~Sternberg}, The infinite groups of
Lie and Cartan, {\it J.~Analyse Math.}~{\bf 15} (1965), 1--114.

\bibitem{Suss} \name{H.~Sussmann,} Orbits of families of vector fields and
integrability of distributions,
{\it Trans.\ Amer.\ Math.~Soc.}~{\bf 180} (1973), 171--188.

\bibitem{Wein} \name{A.~Weinstein}, Symplectic groupoids and Poisson
manifolds, {\it Bull.\ Amer.~Math. Soc.}~{\bf 16}
(1987), 101--104.

\bibitem{28} \name{J. Yorke}, Periods of periodic solutions and the Lipschitz constant, {\it Proc.\ Amer.\ Math.\ Soc.\/}
{\bf 22} (1969), 509--512.

\end{references}
\end{document}